\documentclass[12pt]{article}

\usepackage{amsmath,amssymb,amsfonts}
\usepackage{geometry}
\usepackage[applemac]{inputenc}
\usepackage{amsmath,amssymb}
\usepackage{amsfonts}
\usepackage{amsmath,amssymb}
\usepackage[applemac]{inputenc}
\usepackage{color}
\usepackage{color, colortbl}
\usepackage{xcolor}
\usepackage[most]{tcolorbox}
\usepackage{amsmath}
\usepackage{amssymb}
\usepackage{bbm}
\usepackage{graphicx}
\usepackage{tikz}
\usepackage{geometry}
\usepackage{float}           
\usepackage{graphicx}        
\usepackage{caption}         
\usepackage{subcaption}   
\usepackage{amssymb}
\usetikzlibrary{arrows}
\usepackage{amsfonts}
\usepackage{amsmath,amssymb}
\usepackage[applemac]{inputenc}
\usepackage{color}
\usepackage{amsmath}
\usepackage{amssymb}
\usepackage{graphicx}
\usepackage{tikz}
\usepackage[numbers]{natbib}  
\usepackage{subcaption}
\usepackage{caption}
\newtheorem{theorem}{Theorem}[section]
\newtheorem{lemma}[theorem]{Lemma}
\newcommand{\E}{E}
\newtheorem{proposition}[theorem]{Proposition}
\newtheorem{corollary}[theorem]{Corollary}
\newtheorem{definition}[theorem]{Definition}

\newtheorem{assumption}[theorem]{Assumption}

\newtheorem{problem}[theorem]{Problem}
\newtheorem{remark}[theorem]{Remark}
\usepackage{subcaption}
\usetikzlibrary{arrows.meta}

\usepackage{tikz}
\usetikzlibrary{arrows.meta,patterns,positioning}

\usetikzlibrary{arrows,shapes,automata,petri}

\newcommand{\dproof}{\noindent {Proof.} \quad}
\newcommand{\fproof}{\hfill $\square$ \bigskip}

\numberwithin{equation}{section}

\definecolor{LightCyan}{rgb}{0.88,1,1}

\def\1B{\text{1\!\!I}}

\def\<{\langle}
\def\>{\rangle}

\def\F{\mathcal{F}}

\def\T{\mathcal{T}}

\def\C{\mathcal{C}}
\def\E{\mathbb{E}}
\def\R{\mathbb{R}}
\def\P{\mathbb{P}}

\def\Y{\mathcal{Y}}

\definecolor{strongred}{RGB}{200, 0, 0}

\definecolor{darkpurple}{HTML}{800080}

\definecolor{deepblue}{cmyk}{1,1,0,0.2}

\begin{document}

\begin{center}

{\Large \bf Optimal Stopping for Systems Driven by the Brownian Sheet}

\vspace{0.8cm}

{ by}

\vspace{0.6cm}

\begin{tabular}{cccc}

{ Nacira Agram} &
{ Bernt \O ksendal} &
{ Frank Proske} &
{ Olena Tymoshenko} \\

Royal Institute of Technology &
University of Oslo &
University of Oslo &
University of Oslo \\

Stockholm, Sweden &
Oslo, Norway &
Oslo, Norway &
 NTUU KPI,  Kyiv, Ukraine\\

&&&

\end{tabular}

\end{center}

\vspace{1cm}
\begin{abstract}
We investigate optimal stopping problems for systems driven by the Brownian sheet. 
Our analysis is divided into two parts.\\
In the first part, we derive explicit solutions to two optimal stopping problems 
for the exponentially discounted Brownian sheet. 
The first problem consists in determining the optimal two parameter first hitting point 
$\tau = (\tau_1,\tau_2)$ that maximizes
\[
\mathbb{E}\!\left[e^{-\rho \tau_1 \tau_2} h\big(B(\tau_1,\tau_2)\big)\right],
\]
where $\rho>0$ is a discount factor and $h$ is a reward function. 
Restricting attention to first hitting points of levels, 
we obtain a closed form characterization of the optimal stopping threshold. 
In particular, for linear rewards $h(y)=y$, 
the optimal level is $\widehat y = (2\rho)^{-1/2}$.\\
The second problem concerns optimal stopping of the integrated discounted Brownian sheet:
\[
\mathbb{E}\!\left[\int_0^{\tau_1}\!\!\int_0^{\tau_2}
e^{-\rho tx} B(t,x)\,dt\,dx\right].
\]
We show that the optimal first hitting level is strictly positive,  
and provide an explicit representation of the value function 
in terms of the exponential integral function. 
The optimal threshold is characterized as the unique solution of a nonlinear equation 
derived from a Laplace transform identity for the product $\tau_1\tau_2$.The result is surprising, in view of the (well-known) fact that the  solution of the corresponding one-parameter problem is to stop at the first time $t$ when $B(t)$ hits the negative level $\overline{y}= - (2\rho)^{-1/2}.$\\
In the second and main part of the paper, 
we develop a general potential theoretic framework 
for two parameter optimal stopping problems 
associated with stochastic partial differential equations (SPDEs) driven by the Brownian sheet. 
For a class of It\^o sheets and time space homogeneous SPDEs, 
we prove that the value function is the least superharmonic majorant of the reward 
and that the optimal stopping point is given by the first exit point from the continuation region. 
This yields a general existence theorem for optimal stopping points 
in the plane, extending classical one parameter theory 
to the two parameter Brownian sheet setting.
\end{abstract}



\section{Introduction}

The theory of optimal stopping plays a central role in probability theory and stochastic control, with classical applications to option pricing, sequential analysis, and free boundary problems (see, e.g., Shiryaev~\cite{Shiryaev} and Peskir--Shiryaev~\cite{PeskirShiryaev}). 
Most of the existing literature concerns one parameter processes, in particular Brownian motion and diffusions. 
In contrast, optimal stopping for multi-parameter processes remains far less developed, despite its natural relevance in models involving both time and space.

A fundamental example of a two parameter process is the \emph{Brownian sheet} $B(t,x)$, introduced by Cairoli and Walsh~\cite{CW}. 
It is a centered Gaussian random field with covariance
\[
\mathbb{E}[B(t,x)B(s,y)] = (t \wedge s)(x \wedge y),
\]
exhibiting a spatio-temporal dependence structure that differs substantially from one dimensional Brownian motion. 
Unlike the one parameter case, there is no total ordering of the parameter space $\mathbb{R}_+^2$, which complicates the definition of filtrations, stopping points, martingales, and stochastic integration. 
As a consequence, classical tools such as the optional stopping theorem, Snell envelopes, and excessive majorant characterizations (e.g.\ Dynkin--Yushkevich~\cite{dyk64}) require significant reinterpretation in the two parameter setting.

Two parameter stochastic processes arise naturally in several areas, including image processing, spatial statistics, and financial modeling of random surfaces such as volatility surfaces (see, e.g., Agram et al.~\cite{AOPT5, AOPT4}, Lindgren et al.~\cite{Lin}, and Wackernagel~\cite{WA}). 
In many practical situations, decisions must depend simultaneously on temporal and spatial information. 
For example, in monitoring systems that observe a random environment evolving in both dimensions, determining when and where to stop requires a genuinely two dimensional formulation.

Early developments in two parameter optimal stopping were initiated by Mazziotto and Szpirglas~\cite{mazziotto05}, who constructed a probabilistic framework for stopping rules in $\mathbb{N}^2$ and $\mathbb{R}_+^2$. 
Their approach is used to characterize optimal stopping points via the Snell envelope and to identify them as maximal points up to which the envelope remains a martingale. 
They also considered examples involving the Brownian sheet and bi-Brownian processes, deriving associated systems of variational inequalities. 
Dalang~\cite{Dal} later proved the existence of optimal stopping points for upper semicontinuous two parameter processes using nonstandard probability spaces (Loeb spaces). 
While these works establish existence in abstract settings, they do not provide explicit optimal thresholds nor first exit characterizations in the classical Brownian sheet framework.

More recently, Tanaka~\cite{tanaka90} studied discrete time two parameter stopping and switching problems via dynamic programming, including extensions to Markov settings. 
These contributions provide valuable conceptual foundations but leave open the problem of deriving explicit and analytically tractable solutions in continuous time for the Brownian sheet.

\medskip

{Main contributions.}

This paper develops both explicit solutions and a general potential theoretic framework for optimal stopping problems driven by the Brownian sheet.

\medskip
{(I) Explicit optimal stopping for the discounted Brownian sheet.}

In the first part of the paper we analyze two Brownian sheet analogues of classical stopping problems for Brownian motion.

\begin{enumerate}
\item 
\emph{Optimal stopping of the discounted Brownian sheet.} 
We consider stopping rules given by two parameter stopping points $\tau=(\tau_1,\tau_2)$ corresponding to first hitting points of levels $y$. 
For a discount factor $\rho>0$ and reward function $h$, we study
\[
\Phi(y)=
\mathbb{E}\!\left[
e^{-\rho \tau_1\tau_2}
h\big(B(\tau_1,\tau_2)\big)
\right],
\]
and determine the optimal threshold $\widehat y$ maximizing $\Phi(y)$. 
For linear rewards, we obtain a closed form optimal level and an explicit Laplace transform identity for the product $\tau_1\tau_2$, revealing structural parallels with classical one parameter hitting time theory.

\item 
\emph{Optimal stopping of the integral of the discounted Brownian sheet.} 
We next consider
\[
\Psi(y)=
\mathbb{E}\!\left[
\int_0^{\tau_1(y)}
\int_0^{\tau_2(y)}
e^{-\rho tx} B(t,x)\,dt\,dx
\right],
\]
where $\tau(y)$ denotes the first hitting point of level $y$. We show that the optimal level $y^*$ is strictly positive, 
in contrast to the corresponding one-parameter problem, 
and characterize it as the unique solution of a nonlinear equation 
derived from an explicit Laplace transform representation.
\end{enumerate}

To the best of our knowledge, these are the first closed form optimal stopping thresholds obtained for the discounted Brownian sheet and for its integrated functional.

\medskip

{(II) A potential theoretic framework in the plane.}

In the second and main part of the paper we develop a general theory of optimal stopping for systems driven by the Brownian sheet. 
We consider random fields defined as certain It\^o sheets first and then as solutions of SPDEs of the form
\[
Y(t,x)=
y+
\int_{0}^{t}\int_{0}^{x}
\sigma(Y(u,v))\,B(du,dv),
\]
and formulate the associated two parameter optimal stopping problem in a potential-theoretic framework.

Building on the Cairoli--Walsh theory of stochastic integration on rectangles~\cite{CW}, we introduce suitable notions of superharmonic and excessive functions adapted to the two parameter setting. 
Within this framework we prove a general existence theorem: 
the value function is the least superharmonic majorant of the reward, and the optimal stopping point is characterized as the first exit from the continuation region in the parameter plane.

In contrast to earlier work on optimal stopping for two parameter processes, such as the abstract Snell envelope and compactness-based approaches developed by Giuseppe Mazziotto, Michel Szpirglas ~\cite{mazziotto05}, and Robert C. Dalang~\cite{Dal}, or the theory of bi-Markov processes with commuting semigroups introduced by Mazziotto~\cite{Mazz85}, the optimal stopping point in our setting is identified explicitly as a first exit point from a continuation set in the parameter plane. Existing results for general two parameter processes typically establish existence only via maximality or compactness arguments and do not yield such a first exit characterization. Our approach is genuinely potential theoretic and exploits the specific structure of Brownian sheet driven SDEs, thereby providing a concrete and geometrically transparent alternative to the more abstract techniques prevailing in the literature.

\medskip

{Organization of the paper.}
Section~2 recalls classical one parameter stopping problems for Brownian motion that motivate our extensions. 
Section~3 develops stochastic integration with respect to the Brownian sheet and establishes exponential martingale identities. 
Section~4 solves the first optimal stopping problem for the discounted Brownian sheet. 
Section~5 treats the optimal stopping problem for the integrated discounted sheet. 
Section~6 develops the general potential theoretic theory and proves the existence theorem for optimal stopping of It\^o sheets.

\section{Two classical cases with one parameter Brownian motion}\label{sec.oneparam} 
As mentioned in the Introduction, optimal stopping theory for one parameter stochastic processes has a long history and many applications. As a motivation for our study of optimal stopping for the Brownian sheet, we recall briefly two examples of explicitly solvable classical optimal stopping problems:

\subsection{Optimal stopping of discounted Brownian motion}

Let $B(t)=B(t,\omega); t\ge 0$, $\omega\in\Omega$,
be a one-parameter Brownian motion with $B(0)=0$,
and let $\mathcal T$ denote the set of all stopping times
with respect to the filtration
$\mathbb F=(\mathcal F_t)_{t\ge0}$ generated by $B(\cdot)$.
Let $\rho>0$ be a given discounting constant.

\begin{problem} \label{2.1}
Find a stopping time $\widehat{\tau}$ such that, over all admissible stopping times $\tau\in\mathcal T$,
\begin{align}\label{OST}
   \sup_{\tau \in \mathcal{T}} \E\big[e^{-\rho \tau}B(\tau)\big]
   =
   \E\big[e^{-\rho \widehat{\tau}}B(\widehat{\tau})\big].
\end{align}
\end{problem}
It is well-known (see \cite[Chapter~10]{O})
 that an optimal stopping time $\tau^{*}$ is the \emph{first hitting time} of the form
\begin{align}
    \tau^{*} = \inf\{t\geq 0; B(t)=y^{*}\}
\end{align}
where
\begin{align}
    y^{*}= \frac{1}{\sqrt{2\rho}}.
\end{align}
\begin{remark}The intuition behind this result is that it is natural to wait for a high value of $B(t)$ before stopping, but if one waits too long the discounting will reduce the payoff more than what is gained with high value of $B(t)$. Stopping at the value $y^{*}$ is the optimal tradeoff between these two effects.
\end{remark}

\subsection{Optimal stopping of the integral of the discounted Brownian motion}
A related optimal stopping problem is the following.:
\begin{problem} \label{2.2}
    Find a stopping time $\widehat{\tau}$ such that, over all admissible stopping times $\tau\in\mathcal T$,
\begin{align}\label{OST_int}
   \sup_{\tau \in \mathcal{T}} \E\big[\int_0^{\tau}e^{-\rho t}B(t)dt\big]=\E\big[\int_0^{\widehat{\tau}}e^{-\rho t}B(t)dt\big].
\end{align}
\end{problem}
One can prove (see e.g. \cite{O}, Chapter 10) that an optimal stopping time $\widehat{\tau}$ for this problem is the \emph{first hitting time} of the form
\begin{align}
  \widehat{\tau }= \inf\{t\geq 0; B(t)= \widehat{y}\},
\end{align}
where
\begin{align}
    \widehat{y}= -\frac{1}{\sqrt{2\rho}}. \label{yhat}
\end{align}
\begin{remark} \label{r2.3}
Intuitively, as long as $B(t)>0$ it is optimal to continue, since the integrand in \eqref{OST_int} is positive.  
It is perhaps less obvious that even if $B(t)$ becomes negative, it may still be optimal to wait in the hope that the process returns to positive values later.  
The result above shows that this intuition is correct, but only until $B(t)$ reaches the negative threshold $\widehat{y}$, at which point immediate stopping is optimal.
\end{remark}

These classical one dimensional optimal stopping problems also possess natural financial interpretations, which motivate further generalizations. In the discounted Brownian motion problem \eqref{OST}, for instance, the process $ B(t) $  may represent the price of an asset or the instantaneous return on capital at time $ t $. The discount factor $ e^{-\rho \tau} $ reflects the time value of money. In this case, the optimal stopping time corresponds to deciding when it is most beneficial to stop and lock in the current gain, since waiting too long reduces the payoff due to discounting.

Similarly, the stopping problem \eqref{OST_int} aims at maximizing the total expected discounted return accumulated over time. Here, the decision to stop depends on whether the expected continuation value still compensates for the increasing discount penalty.
 
Financial interpretations of one dimensional optimal stopping problems are classical (see \cite{Marcozzi01}, \cite{PeskirShiryaev}), while two parameter models can describe more complex investment strategies involving multiple funding streams or sources \cite{Lin}.
In our setting, we consider a two parameter generalization where the underlying process is a Brownian sheet $B(t, x) $, and the objective is to maximize the expected value
 $$\E[e^{-\rho \tau_1 \tau_2} h(B(\tau_1,\tau_2))],$$
or  
\[
\mathbb{E}\left[\int_0^{\tau_1} \int_0^{\tau_2} e^{-\rho s x} B(s, x) \, ds \, dx\right],
\]
where $ (\tau_1, \tau_2) $ is a stopping point in the two dimensional domain and $h$ is a reward function.

This double integral can be understood as the total accumulated discounted value of a cash flow that depends on two variables. The first variable $ t $ represents time. The second variable $ x $ represents the flow of funds, such as the volume of investments entering the system or cash inflows from external sources.

The exponential weight $ e^{-\rho s x} $ applies discounting with respect to both time and flow magnitude. This reflects that later or larger flows have less present value.

Such models may arise in finance when value depends jointly on time and the size or source of incoming capital. Examples include investment decisions involving multiple funding streams or managing distributed financial inflows over time and scale.

\section{Exponential of Brownian sheet integrals}\label{sec:exp-brownian-sheet}
To study optimal stopping problems in the two parameter setting with the Brownian sheet, we first need some key results from stochastic calculus in the plane, which provides a powerful extension of classical It\^o  calculus to systems driven by two parameter noise.To handle stochastic integration and optimal stopping in the two parameter setting, we rely on the theory of stochastic integrals in the plane, developed by Cairoli and Walsh \cite{CW}, and related foundational works \cite{WZ}, \cite{Wo}.
The Brownian sheet, which generalizes standard Brownian motion to two spatial directions, exhibits complex dependency structures and requires a more elaborate integration theory. 

In this section, we provide a detailed derivation of the exponential formula for functionals of Brownian sheet integrals, including a rigorous proof of Proposition~\ref{cor1}. Our formulation highlights the structure and causal interaction encoded by the indicator $I(\zeta \bar{\wedge} \zeta')$, offering a clear and accessible expression of the two parameter exponential expansion. This representation serves as a building block for later developments in optimal stopping and stochastic control in space time domains.

We denote by $\{B(t,x);\ t,x \geq 0 \}$ a Brownian sheet and $(\Omega, \mathcal{F}, P)$ a complete probability space on which we define the (completed) filtration $\{\mathcal{F}_{t,x}\}$ generated by $B(s,a)$, $s \leq t$, $a \leq x$. 

All stochastic integrals with respect to the Brownian sheet are understood in the sense of Cairoli \& Walsh integration theory (see \cite{CW}).
 
Let us consider the stochastic field
\begin{equation}\label{Y_def}
Y(z) = \int_{R_z} \alpha(\zeta)\,d\zeta + \int_{R_z} \beta(\zeta)\,B(d\zeta), \quad z = (t,x) \in \mathbb{R}^2_+,
\end{equation}
where $R_z = [0,t] \times [0,x]$ and $\alpha, \beta: \mathbb{R}^2_+ \to \mathbb{R}$ are bounded, $\F_z$-adapted processes. The notation 
$B(d\zeta)$ stands for the stochastic integral with respect to the Brownian sheet over a differential element 
$d\zeta=dsdx$, and the integration is taken over the rectangle $R_z$. Define the exponential process
\begin{equation}\label{Gz_def}
    G(z) := e^{Y(z)} = \exp\left( \int_{R_z} \alpha(\zeta)\,d\zeta + \int_{R_z} \beta(\zeta)\,B(d\zeta) \right).
\end{equation}
Before proceeding, define the indicator function $I$ as follows:
\begin{equation}\label{ind}
    I((a_1, a_2) \bar{\wedge} (b_1, b_2)) =
\begin{cases} 
1, & \text{if } a_1 \leq b_1 \text{ and } a_2 \geq b_2, \\
0, & \text{otherwise}.
\end{cases}
\end{equation}

To make sure that the integrals in $Y(z)$ and $G(z)$ are well defined, we require the following assumptions on $\alpha$ and $\beta$.
\begin{assumption}\label{a_int}
We assume that the coefficient fields $\alpha,\beta:\R^2_+\to\R$ are jointly measurable, predictable 
with respect to the filtration $\{\F_{t,x}\}$, and satisfy for every finite rectangle $R=[0,t]\times[0,x]$:
\[
\E\int_R |\alpha(\zeta)|\,d\zeta < \infty, 
\qquad 
\E\int_R \beta(\zeta)^2\,d\zeta < \infty.
\]
In addition, we assume
\[
\E \exp\Big(\tfrac12\int_R \beta(\zeta)^2\,d\zeta\Big)<\infty,
\]
ensuring integrability of the exponential process $G(z)=e^{Y(z)}$.
\end{assumption}

We now establish the two parameter analogue of the exponential martingale formula,  valid for the stochastic exponential driven by the Brownian sheet.

\begin{proposition}[two parameter Exponential Formula]\label{cor1}
Let $Y(z)$ and $G(z)$ be the stochastic exponential fields defined by \eqref{Y_def} and \eqref{Gz_def} respectively. Then  $G(z)$ admits the following  representation
\begin{align*}
G(z) &= 1 + \int_{R_z} G(\zeta)\left( \alpha(\zeta) + \tfrac{1}{2} \beta^2(\zeta) \right) d\zeta \\
&\quad + \iint_{R_z \times R_z} I(\zeta \bar{\wedge} \zeta')\, G(\zeta \vee \zeta')\, \left( \alpha(\zeta) + \tfrac{1}{2} \beta^2(\zeta) \right)\left( \alpha(\zeta') + \tfrac{1}{2} \beta^2(\zeta') \right) d\zeta\, d\zeta' \\
&\quad + \text{terms with zero expectation}.
\end{align*}
\end{proposition}

 To derive an explicit expansion for the exponential process  $G(z)$ defined in \eqref{Gz_def}, we apply the two parameter It\^o  formula as developed in \cite{AOPT1}, see also \cite{CW}.

\subsection{The two parameter It\^o formula}

We recall the two parameter It\^o formula following Wong and Zakai~\cite{WZ}.
Throughout, we use the following notation.

\begin{figure}[t]
\centering

\begin{subfigure}[t]{0.32\textwidth}
\centering
\begin{tikzpicture}[scale=0.70, >=Latex]
  \def\tzero{1.2}
  \def\tend{5.0}
  \def\xzero{1.0}
  \def\xend{3.2}

  \colorlet{rectfill}{blue!15}
  \colorlet{rectborder}{blue!60!black}
  \colorlet{cornercol}{blue!80!black}

  \draw[->] (0,0) -- (6,0);
  \draw[->] (0,0) -- (0,4.2);

  \draw (\tzero,0) -- (\tzero,-0.1);
  \draw (\tend,0) -- (\tend,-0.1);
  \draw (0,\xzero) -- (-0.1,\xzero);
  \draw (0,\xend) -- (-0.1,\xend);

  \node[below] at (\tzero,0) {$t_0$};
  \node[below] at (\tend,0) {$t$};
  \node[left]  at (0,\xzero) {$x_0$};
  \node[left]  at (0,\xend) {$x$};

  \fill[rectfill] (\tzero,\xzero) rectangle (\tend,\xend);
  \draw[thick,rectborder] (\tzero,\xzero) rectangle (\tend,\xend);

  \fill[cornercol] (\tzero,\xzero) circle (2.5pt);
  \fill[cornercol] (\tend,\xend) circle (2.5pt);

  \node[below left,cornercol] at (\tzero,\xzero) {$z_0$};
  \node[above right,cornercol] at (\tend,\xend) {$z$};

  \node[rectborder] at (3.1,2.1) {$R_{z_0}(z)$};

\end{tikzpicture}
\caption{Rectangle $R_{z_0}(z)$.}
\label{fig:rectangle}
\end{subfigure}\hfill


\begin{subfigure}[t]{0.34\textwidth}
\centering
\begin{tikzpicture}[scale=0.70, >=Latex]
  \draw[->] (0,0) -- (6,0);
  \draw[->] (0,0) -- (0,4.2);

  \def\s{1.5}
  \def\sp{4.5}
  \def\a{3.0}
  \def\ap{1.2}

  \draw (\s,0)  -- (\s,-0.1);
  \draw (\sp,0) -- (\sp,-0.1);
  \draw (0,\a)  -- (-0.1,\a);
  \draw (0,\ap) -- (-0.1,\ap);

  \node[below] at (\s,0) {$s$};
  \node[below] at (\sp,0) {$s'$};
  \node[left]  at (0,\a) {$a$};
  \node[left]  at (0,\ap) {$a'$};

  \colorlet{zetacol}{blue!70!black}
  \colorlet{zetapcol}{red!70!black}
  \colorlet{maxcol}{purple!80!black}

  \fill[zetacol] (\s,\a) circle (2pt);
  \fill[zetapcol] (\sp,\ap) circle (2pt);

  \node[above left,xshift=28pt, yshift=1pt, zetacol] at (\s,\a) {$\zeta=(s,a)$};
  \node[below right,zetapcol] at (\sp,\ap) {$\zeta'=(s',a')$};

  \fill[maxcol] (\sp,\a) circle (2pt);
  \node[above right,maxcol] at (\sp,\a) {$\zeta\vee\zeta'$};

  \draw[dashed,maxcol] (\s,\a) -- (\sp,\a);
  \draw[dashed,maxcol] (\sp,\ap) -- (\sp,\a);
\end{tikzpicture}
\caption{Geometry of $\zeta\vee\zeta'$.}
\label{fig:zeta-vee}
\end{subfigure}\hfill
\begin{subfigure}[t]{0.32\textwidth}
\centering
\begin{tikzpicture}[scale=0.85, >=Latex]

  \def\xmax{5.5}
  \def\ymax{4.2}

  \def\bOne{3.2}
  \def\bTwo{2.2}

  \draw[->] (0,0) -- (\xmax,0) node[right] {$a_1$};
  \draw[->] (0,0) -- (0,\ymax) node[above] {$a_2$};

  \fill[blue!18] (0,\bTwo) rectangle (\bOne,\ymax);

  \draw[dashed] (\bOne,0) -- (\bOne,\ymax);
  \draw[dashed] (0,\bTwo) -- (\xmax,\bTwo);

  \fill[red!80!black] (\bOne,\bTwo) circle (2.5pt);
  \node[above right] at (\bOne,\bTwo)
  {$b$};

  \node[left] at (2.5,3.6) {\small $a_1 \le b_1$};
  \node[left] at (2.5,3.1) {\small $a_2 \ge b_2$};

\end{tikzpicture}
\caption{Region where $I=1$.}
\label{fig:mixed-order}
\end{subfigure}

\caption{
Geometric interpretation of rectangle integration,
the coordinatewise maximum $\zeta\vee\zeta'$,
and the mixed partial order.
}
\label{fig:two-figs}
\end{figure}
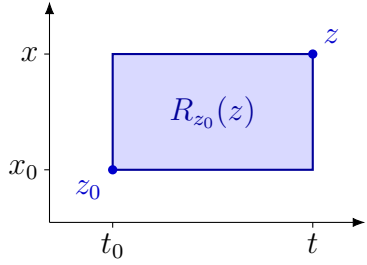
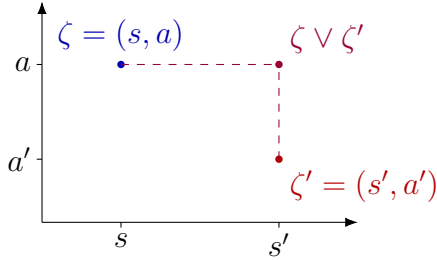
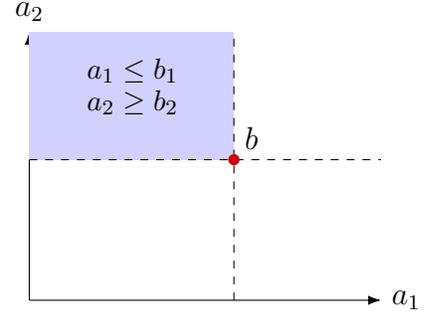

\begin{itemize}

\item
We write $\zeta=(\zeta_1,\zeta_2)=(s,a)\in \mathbb{R}^2$ and
\[
d\zeta = d\zeta_1\, d\zeta_2 = ds\,da.
\]

\item
$B(t,x)$ denotes a Brownian sheet, $t\ge0$, $x\in\mathbb{R}$.

\item
For $z=(z_1,z_2)=(t,x)$ we define the rectangle
\[
R_z=[0,z_1]\times[0,z_2],
\]
see Figure~\ref{fig:rectangle}.

\item
For an adapted process $\varphi$,
\[
\int_{R_z}\varphi(\zeta)\,B(d\zeta)
\]
denotes the It\^o integral with respect to $B$ over $R_z$.

\item
For a measurable function $\psi$,
\[
\int_{R_z}\psi(\zeta)\,d\zeta
\]
denotes the two dimensional Lebesgue integral over $R_z$.

\item
For $a=(a_1,a_2)$ and $b=(b_1,b_2)$ we define
\[
a\vee b
=
(\max\{a_1,b_1\},\,\max\{a_2,b_2\}),
\]
whose geometric interpretation is illustrated in
Figure~\ref{fig:zeta-vee}.

\item
We introduce the mixed partial order
\[
I((a_1,a_2)\bar{\wedge}(b_1,b_2))
=
\mathbbm{1}_{\{a_1 \le b_1,\; a_2 \ge b_2\}},
\]
whose geometry is shown in Figure~\ref{fig:mixed-order}.

\end{itemize}

\begin{theorem}[It\^o formula]\label{Ito}
Suppose that the process $Y(z)=Y(t,x)$, $z=(t,x)$, satisfies the SPDE
\begin{equation}\label{Y}
  dY(z)=\alpha(z)\,dz+\beta(z)\,B(dz),
\end{equation}
with given initial values
\begin{equation*}\label{Y-ic}
  Y(t_0,x_0)=y,\qquad D_1Y(t,x_0)=p(t,x_0)\ (t>t_0),\qquad D_2Y(t_0,x)=q(t_0,x)\ (x>x_0).
\end{equation*}
Equivalently,
\begin{equation}\label{Y-pde}
  \frac{\partial^2}{\partial t\,\partial x}Y(t,x)=\alpha(t,x)+\beta(t,x)\diamond W(t,x),
\end{equation}
where $\alpha$ and $\beta$ are integrable adapted processes, $\diamond$ denotes the Wick product, and
$W(t,x)=\frac{\partial^2}{\partial t\,\partial x}B(t,x)$ is the space time white noise. Moreover, for (\ref{Y-pde})
\[
  D_1Y(t,x_0)=\frac{\partial}{\partial t}Y(t,x_0),\qquad
  D_2Y(t_0,x)=\frac{\partial}{\partial x}Y(t_0,x).
\]
In integrated form,
\begin{align}\label{Y-int}
Y(t,x)=\;& y+\int_{t_0}^{t} p(s,x_0)\,ds+\int_{x_0}^{x} q(t_0,a)\,da \nonumber\\
&+\int_{t_0}^{t}\int_{x_0}^{x}\alpha(s,a)\,ds\,da
+\int_{t_0}^{t}\int_{x_0}^{x}\beta(s,a)\,B(ds\,da),
\qquad t\ge t_0,\ x\ge x_0 .
\end{align}
Let $f:\R \mapsto \R$ be in $\C^4(\mathbb{R})$, i.e. $f$ is four times continuously differentiable.\\ 
Then, with $z=(t,x)$ and $z_0=(t_0,x_0)$,
\begin{align*}
f(Y(z)) &= f(Y(z_0)) +\int_{t_0}^t f'(Y(s,x_0))D_1 Y(s,x_0)ds
+ \int_{x_0}^x f'(Y(t_0,a)) D_2 Y(t_0,a)da \nonumber\\
&+ \int_{R_{z_0}(z)} f'(Y(\zeta))\left[\alpha(\zeta)\,d\zeta + \beta(\zeta)\,B(d\zeta)\right] \\
&\quad + \frac{1}{2} \int_{R_{z_0}(z)} f''(Y(\zeta))\, \beta^2(\zeta)\, d\zeta \\
&\quad + \iint_{R_{z_0}(z) \times R_{z_0}(z)}f''(Y(\zeta \vee \zeta'))\, \beta(\zeta)\beta(\zeta')\, B(d\zeta)\, B(d\zeta') \\
&\quad + \iint_{R_{z_0}(z) \times R_{z_0}(z)} f''(Y(\zeta \vee \zeta'))\, \beta(\zeta')\alpha(\zeta)\, d\zeta\, B(d\zeta') \\
&\quad + \iint_{R_{z_0}(z) \times R_{z_0}(z)} f''(Y(\zeta \vee \zeta'))\, \beta(\zeta)\alpha(\zeta')\, B(d\zeta)\, d\zeta' \\
&\quad + \iint_{R_{z_0}(z) \times R_{z_0}(z)} I(\zeta \bar{\wedge} \zeta')\, f''(Y(\zeta \vee \zeta'))\, \alpha(\zeta)\alpha(\zeta')\, d\zeta\, d\zeta' \\
&\quad + \frac{1}{2} \iint_{R_{z_0}(z) \times R_{z_0}(z)} I(\zeta \bar{\wedge} \zeta')\, f^{(3)}(Y(\zeta \vee \zeta')) \left[ \alpha(\zeta)\beta^2(\zeta') + \alpha(\zeta')\beta^2(\zeta) \right]\, d\zeta\, d\zeta' \\
&\quad + \frac{1}{4} \iint_{R_{z_0}(z) \times R_{z_0}(z)} I(\zeta \bar{\wedge} \zeta')\, f^{(4)}(Y(\zeta \vee \zeta'))\, \beta^2(\zeta)\, \beta^2(\zeta')\, d\zeta\, d\zeta',
\end{align*}
where $R_{z_0}(z)=R_z \setminus R_{t_0,x} \setminus R_{t,z_0} \cup R_{t_0,x_0}$; $\zeta \vee \zeta' = \bigl(\max\{s,s'\},\, \max\{a,a'\}\bigr)$,  and the indicator function $I$ is defined in \eqref{ind}.
\end{theorem}
In the particular case where $f(y) = e^y$ all derivatives satisfy $f^{(k)}(y) = e^y$ for every $k \geq 0$. Therefore, if we in addition assume that $ t_0=x_0=0=Y(0,x)=Y(t,0); \  \forall t > 0, x >0$,the expansion simplifies significantly, as all derivative terms reduce to powers of the same exponential function. Since $G(z) = e^{Y(z)}$, with $Y(z)$ defined by~\eqref{Y_def}, we have
\[
f'(Y(\zeta)) = f''(Y(\zeta)) = f^{(3)}(Y(\zeta)) = f^{(4)}(Y(\zeta)) = G(\zeta).
\]
Substituting into the general expansion yields:
\begin{align*}
G(z) &= G(0) 
+ \int_{R_z} G(\zeta)\, \alpha(\zeta)\, d\zeta 
+ \frac{1}{2} \int_{R_z} G(\zeta)\, \beta^2(\zeta)\, d\zeta \\
&\quad + \iint_{R_z \times R_z} G(\zeta \vee \zeta')\, \beta(\zeta)\beta(\zeta')\, B(d\zeta)\, B(d\zeta') \\
&\quad + \iint_{R_z \times R_z} G(\zeta \vee \zeta')\, \beta(\zeta')\alpha(\zeta)\, d\zeta\, B(d\zeta') \\
&\quad + \iint_{R_z \times R_z} G(\zeta \vee \zeta')\, \beta(\zeta)\alpha(\zeta')\, B(d\zeta)\, d\zeta' \\
&\quad + \iint_{R_z \times R_z} I(\zeta \bar{\wedge} \zeta')\, G(\zeta \vee \zeta')\, \alpha(\zeta)\alpha(\zeta')\, d\zeta\, d\zeta' \\
&\quad + \frac{1}{2} \iint_{R_z \times R_z} I(\zeta \bar{\wedge} \zeta')\, G(\zeta \vee \zeta')\left[\alpha(\zeta)\beta^2(\zeta') + \alpha(\zeta')\beta^2(\zeta)\right]\, d\zeta\, d\zeta' \\
&\quad + \frac{1}{4} \iint_{R_z \times R_z} I(\zeta \bar{\wedge} \zeta')\, G(\zeta \vee \zeta')\, \beta^2(\zeta)\beta^2(\zeta')\, d\zeta\, d\zeta'.
\end{align*} 
It was shown in \cite{WZ} that all stochastic integrals involving the Brownian sheet, including the double, triple, and quadruple integrals, as well as the mixed integrals $B(d\zeta)\,d\zeta'$ and $d\zeta\,B(d\zeta')$, are weak martingales. Consequently, under suitable integrability and adaptedness conditions, their expectations vanish.

Therefore, taking expectations and using the zero-mean property of all stochastic terms, we obtain:
\begin{align*}
\mathbb{E}[G(z)] &= G(0) 
+ \mathbb{E}\left[ \int_{R_z} G(\zeta) \alpha(\zeta)\,d\zeta + \frac{1}{2} \int_{R_z} G(\zeta) \beta^2(\zeta)\,d\zeta \right] \\
& + \mathbb{E}\left[ \iint_{R_z \times R_z} I(\zeta \bar{\wedge} \zeta')\, G(\zeta \vee \zeta')\, \left( \alpha(\zeta) + \tfrac{1}{2} \beta^2(\zeta) \right)\left( \alpha(\zeta') + \tfrac{1}{2} \beta^2(\zeta') \right)\, d\zeta\,d\zeta' \right].
\end{align*}


In particular, choosing $\alpha= -\tfrac{1}{2} \beta^2$ we get
\begin{proposition}\label{cor2}
Let $\beta$ be a bounded and $\F_z$-adapted process. Then
\[
\mathbb{E}\left[ \exp\left( \int_{R_{z}}  \beta(\zeta) B(d\zeta) - \tfrac{1}{2}\int_{R_{z}}\beta^2(\zeta) d\zeta \right) \right] = 1, \quad z=(t,x) \in \mathbb{R}_{+}^2.
\]
This identity generalizes the classical exponential martingale property to the Brownian sheet setting.
\end{proposition}

\begin{remark}
The stochastic integral in Proposition~\ref{cor2}, namely
\[
\int_{R_z} \beta(\zeta)\, B(d\zeta),
\]
is interpreted in the It\^ o sense as defined in the current paper. This formulation differs from the Walsh-type integral developed in \cite{Walsh1986}, or the Cairoli-Walsh framework and related approaches in \cite{DalangMueller2009} and \cite{DalangBook}. In particular, our setting focuses on exponential martingales over rectangular domains in \( \mathbb{R}_+^2 \), allowing us to derive a Dol\' eans-Dade-type exponential in the two parameter setting that generalizes the classical one dimensional result.
\end{remark}

\section{The First Optimal Stopping Problem: Discounted Brownian sheet}\label{sec:first-optimal-stopping}

We now proceed to study the first optimal stopping problem. We first define the notion of stopping points in the two parameter setting.

\begin{definition}
Let $\mathcal{F}_{t,x}$ denote the $\sigma$-algebra generated by $$\{B(r,a); \ 0\leq r \leq t, 0 \leq a \leq x\}.$$ A random variable $\tau = (\tau_1(\omega), \tau_2(\omega)) \colon \Omega \to \mathbb{R}_+ \times \mathbb{R}_+$ is called a \emph{stopping point} if
\[
\{\omega \in \Omega : \tau_1(\omega) \leq t \ \& \ \tau_2(\omega) \leq x\} \in \mathcal{F}_{t,x}, \quad \forall (t,x) \in  \mathbb{R}_+^2.
\]
This means that it is possible to determine whether the stopping has occurred by observing the Brownian sheet up to $(t,x)$.

We denote the set of all such stopping points by $\mathcal{T}$.\\
Note that it is possible to have $\tau_i(\omega)=\infty$ for some (or all) $\omega \in \Omega$; i=1,2.
\end{definition}
Next, we show that the stochastic integral from Corollary~\ref{cor2} remains well defined when the domain of integration is randomized via a stopping point.
\begin{proposition}
 Let $\varphi \colon \mathbb{R}_+^2 \times \Omega \to \mathbb{R}$ be a bounded predictable process with respect to the filtration $\{\mathcal{F}_{t,x}\}$. If $\tau = (\tau_1, \tau_2)$ is a stopping point such that $\mathbb{E}[\tau_1 \tau_2] < \infty$, then the stochastic integral
\[
\int_0^{\tau_1} \int_0^{\tau_2} \varphi(\zeta)\, B(d\zeta)
\]
is well defined in \(L^2(\Omega)\), and it satisfies the isometry
\[
\mathbb{E} \left[ \left( \int_0^{\tau_1} \int_0^{\tau_2} \varphi(\zeta)\, B(d\zeta) \right)^2 \right] = 
\mathbb{E} \left[ \int_0^{\tau_1} \int_0^{\tau_2} \varphi^2(\zeta)\, d\zeta \right] < \infty.
\]
\end{proposition}

\dproof
Since the process $\varphi$ is bounded and adapted with respect to the filtration $\{\mathcal{F}_{t,x}\}$, the stochastic integral
\[
\int_0^{\tau_1} \int_0^{\tau_2} \varphi(\zeta)\, B(d\zeta)
\]
can be written as an integral over the random rectangle $[0,\tau_1] \times [0,\tau_2]$. More precisely,
\[
\int_0^{\tau_1} \int_0^{\tau_2} \varphi(\zeta)\, B(d\zeta) = \int_{\mathbb{R}_+^2} \mathbbm{1}_{\{ \zeta \le \tau \}}\, \varphi(\zeta)\, B(d\zeta),
\]
where $\zeta \le \tau$ means $\zeta_1 \le \tau_1$ and $\zeta_2 \le \tau_2$.
To compute the second moment, we apply the It\^o isometry for the Brownian sheet
\begin{align}
\mathbb{E} \left[ \left( \int_0^{\tau_1} \int_0^{\tau_2} \varphi(\zeta)\, B(d\zeta) \right)^2 \right]
&= \mathbb{E} \left[ \int_0^{\tau_1} \int_0^{\tau_2} \varphi^2(\zeta)\, d\zeta \right]. \label{L2_res}
\end{align}
To justify this identity rigorously, we consider truncated versions: for $T, X > 0$, define $\tau_1^T = \tau_1 \wedge T$, $\tau_2^X = \tau_2 \wedge X$, and use
\[
\mathbb{E} \left[ \left( \int_0^{\tau_1^T} \int_0^{\tau_2^X} \varphi(\zeta)\, B(d\zeta) \right)^2 \right]
= \mathbb{E} \left[ \int_0^{\tau_1^T} \int_0^{\tau_2^X} \varphi^2(\zeta)\, d\zeta \right].
\]
Since $\varphi$ is bounded and the integrand increases pointwise as $T,X \to \infty$, we apply the \emph{monotone convergence theorem} (Beppo Levi) to obtain
\[
\mathbb{E} \left[ \int_0^{\tau_1} \int_0^{\tau_2} \varphi^2(\zeta)\, d\zeta \right] = \lim_{T,X \to \infty}
\mathbb{E} \left[ \int_0^{\tau_1^T} \int_0^{\tau_2^X} \varphi^2(\zeta)\, d\zeta \right] < \infty.
\]
Thus, the stochastic integral is well-defined in $L^2(\Omega)$ and satisfies~\eqref{L2_res}. Moreover, the integral has zero expectation 
\[
\mathbb{E} \left[ \int_0^{\tau_1} \int_0^{\tau_2} \varphi(\zeta)\, B(d\zeta) \right] = 0.
\]
\fproof 

\begin{definition}
A stopping point $\tau \in \mathcal{T}$ is called a \emph{first hitting point for the level $y \in \mathbb{R}$} if
\[
B(\tau) = y \quad \text{a.s., and} \quad B(z) \neq y \ \text{for all } z \leq \tau, \ z \neq \tau.
\]
We denote the set of all such stopping points by $\mathcal{T}_y$. Note that first hitting points need not be unique.

We let 
    $\widehat{\T}=\cup_{y\in \R} \T_y$
denote the set of all first hitting points.
\end{definition}

\begin{remark}
Let us consider $Y(z) := B(z)$ and $f(y) := y^2$. Then the It\^o formula in two dimensions yields
\[
\mathbb{E}[B^2(t,x)] = \int_{[0,t]\times[0,x]} d\zeta = t x.
\]
Fix $R > 0$, Assume that there exists a first hitting point $\tau^{(R)} \in T_R \cup T_{-R}$, and fix one such choice. Then
 $|B(\tau^{(R)})| = R$ almost surely (a.s.), and
\[
\mathbb{E}[\tau^{(R)}_1 \tau^{(R)}_2] = R^2 < \infty.
\]
Thus, the It\^o formula remains valid for stopping points $\tau$ such that $\tau \leq \tau^{(R)}$ a.s. for some $R$.
\end{remark}

We begin with the following auxiliary result:       
\begin{proposition}\label{lemma1}
Let $y \in \R$, $\beta >0$ and let $\tau_y=(\tau_{y,1},\tau_{y,2})$ be the first hitting point of the level $y$ (i.e.\ $B(\tau_y)=y$ a.s.). Then
    \begin{align}
\E\left[ e^{-\tfrac{1}{2} \beta^2 \tau_{y,1}\tau_{y,2}} \right] =  e^{-\beta|y|}. \label{tau}
\end{align}
\end{proposition}
\dproof
(a) First assume that $y>0$. 
Choose $R>0$ and define $\tau^{(R)}$ to be the first exit point of 
$B(t,x)$ from the interval $[-R,y]$.
Since $B(\tau^{(R)}) \in \{-R,y\}$, we have
\[
0 = \mathbb{E}[B(\tau^{(R)})]
= -R \mathbb{P}(B(\tau^{(R)})=-R)
+ y \mathbb{P}(B(\tau^{(R)})=y).
\]
Hence
\[
R\,\mathbb{P}(B(\tau^{(R)})=-R) \le y,
\]
and therefore 
$\mathbb{P}(B(\tau^{(R)})=-R)\to 0$ as $R\to\infty$.
Thus $\tau_y<\infty$ a.s. and $B(\tau_y)=y$.
By Proposition \ref{cor2}  applied to the stopped exponential martingale,
\[
1 =
\mathbb{E}\!\left[
\exp\!\left(
\beta B(\tau^{(R)})
-\tfrac12 \beta^2 \tau^{(R)}_1\tau^{(R)}_2
\right)
\right].
\]
Splitting according to the exit location,
\begin{align*}
1
&= e^{\beta y}
\mathbb{E}\!\left[
e^{-\frac12 \beta^2 \tau^{(R)}_1\tau^{(R)}_2}
\mathbbm{1}_{\{B(\tau^{(R)})=y\}}
\right]  \\
&\quad + e^{-\beta R}
\mathbb{E}\!\left[
e^{-\frac12 \beta^2 \tau^{(R)}_1\tau^{(R)}_2}
\mathbbm{1}_{\{B(\tau^{(R)})=-R\}}
\right].
\end{align*}
Since $0 \le e^{-\beta R}\le 1$ and $e^{-\beta R}\to 0$, 
the second term vanishes as $R\to\infty$.
Passing to the limit yields
\[
1 = e^{\beta y}
\mathbb{E}\!\left[
e^{-\frac12 \beta^2 \tau_{y,1}\tau_{y,2}}
\right].
\]
Therefore,
\[
\mathbb{E}\!\left[
e^{-\frac12 \beta^2 \tau_{y,1}\tau_{y,2}}
\right]
= e^{-\beta y}.
\]
(b) If $y<0$, the result follows by symmetry of the Brownian sheet.\\
(c) If $y=0$, the identity is immediate.
\fproof
\begin{remark}
The Proposition above provides an explicit Laplace transform of the product
$\tau_{y,1}\tau_{y,2}$ of the first hitting coordinates. 
To our knowledge, such closed form identities are rare in
two parameter stopping problems.\end{remark}

\subsection*{Solution of the first problem: Optimal stopping of the discounted Brownian sheet}
We  consider the optimal stopping problem formulated in terms of first hitting points. Given a reward function 
$h:\mathbb{R}_+ \to \mathbb{R}_+$ and a discount rate $\rho>0$, the goal is to maximize the expected discounted reward. This leads to the definition of the value function 
$\Phi$, characterized by the following result:
\begin{theorem} [Characterization of optimal first hitting point]
Let $h:\mathbb{R}_+ \to \mathbb{R}_{+}$ be a $C^1$ function and let $\rho > 0$ be a discounting constant.
Define the function $g:\R_+\to\R_+$ by
\begin{align}
g(y)= h(y) \cdot e^{-\sqrt{2\rho} y}. \label{max}
\end{align}
\begin{enumerate}
\item
 Define the first hitting value function by
\begin{align}
\Phi(y) := \sup_{\tau \in \mathcal{T}_y} \mathbb{E}\left[ e^{-\rho \tau_1 \tau_2} h\left(B(\tau_1,\tau_2)\right) \right], \quad y>0.
\end{align}
Then $\Phi(y)=g(y)$. \\
Hence $\Phi(y)$ attains its maximum at $\hat{y}$ if and only if $y=\hat{y}$ is a maximum point of $g(y)$.

Therefore, if $\Phi$ has a maximum point $\hat{y}$, then $\hat{y}$ satisfies the equation
\begin{align}
h'(\hat{y}) = \sqrt{2\rho} \, h(\hat{y}).
\end{align}

\item
In particular, for the choice $h(y) = y$, we obtain: \\
The solution to the optimal stopping problem
\begin{align}
\Phi_1(y) := \sup_{\tau \in \mathcal{T}_y} \mathbb{E}\left[ e^{-\rho \tau_1 \tau_2} B(\tau_1, \tau_2) \right], \quad y >0,
\end{align}
is attained at
\begin{align}
\hat{y} = \frac{1}{\sqrt{2\rho}}.
\end{align}
This yields the maximum value
\begin{align}
\sup_{y \in \mathbb{R}_+} \Phi_1(y) = \Phi_1(\hat{y}) = \frac{1}{\sqrt{2\rho}} e^{-1}. \label{sup}
\end{align}
\end{enumerate}
\end{theorem}
\dproof Fix $y>0$ and let $\tau=(\tau_1,\tau_2)\in \mathcal{T}_y$. Since $B(\tau_1,\tau_2)=y$ a. s., we have
\[
h\!\big(B(\tau_1,\tau_2)\big)=h(y)\quad\text{a.s.}
\]
Hence,
\begin{align*}
\E[e^{-\rho \tau_1\tau_2} h(B(\tau_1,\tau_2))]= h(y) \cdot \E[e^{-\rho \tau_1\tau_2}].
\end{align*}
Using Proposition \ref{lemma1} with $\beta = \sqrt{2\rho}$ yields
\begin{align}
\E[e^{-\rho \tau_1\tau_2} h(B(\tau_1,\tau_2))]= h(y) \cdot e^{-\sqrt{2\rho} y} = g(y). \label{max}
\end{align}
If $\Phi$ attains its maximum at an interior point $\hat y>0$, then 
\begin{align*}
0=g'(y) = h'(y)e^{-\sqrt{2\rho} y} - \sqrt{2\rho} h(y)e^{-\sqrt{2\rho} y}.
\end{align*}
We conclude that if there is an optimal first hitting level $y$ then
\begin{align}
h'(\hat{y}) = \sqrt{2\rho} \, h(\hat{y}). \label{haty}
\end{align}
If $h(y)=y$, then \eqref{haty} gives the unique solution $\hat{y} = \frac{1}{\sqrt{2\rho}}$, which substituted into \eqref{max} yields the maximum value in \eqref{sup}.

\fproof

\begin{remark}
    Note that equation \eqref{haty} does not have a solution $\hat{y}$ for all functions $h$. For example, if $h(y)= e^{ky}$ then there is no solution. Indeed, if $h(y)=e^{ky}$, then
    $g(y)=\exp(ky-\sqrt{2\rho} y)$
    which has no maximum point in $(0,\infty)$.
\end{remark}

\begin{corollary}[Linear Reward]
Suppose \(h(y) = y\). Then the optimal stopping threshold is \(\hat{y} = \frac{1}{\sqrt{2\rho}}\), and the maximal expected value is
\[
\Phi_1(\hat{y}) = \frac{1}{\sqrt{2\rho}} e^{-1}.
\]
\end{corollary}

\begin{remark}[Comparison with one parameter Case]
In the classical one parameter Brownian motion case, the optimal stopping time for a problem of the form
\[
\sup_{\tau} \mathbb{E}[e^{-\rho \tau} h(B(\tau))]
\]
typically yields the condition \(h'(y) = \sqrt{2\rho} \, h(y)\) as well. However, the discounting is exponential in \(\tau\), while in the two parameter case it is exponential in the product \(\tau_1 \tau_2\), leading to a more delicate balance between the spatial and temporal components. Moreover, the Brownian sheet lacks the strong Markov property in the same form, and first hitting points are not uniquely defined making the analysis more subtle.
\end{remark}
\begin{corollary}
Let $h(y) = y^n$ for some integer $n \in \mathbb{N}$, $n \geq 1$. Then the optimal first hitting level $\hat{y}$ that maximizes the expected discounted reward
\[
\Phi_n(y) := \sup_{\tau \in \mathcal{T}_y} \mathbb{E}\left[ e^{-\rho \tau_1 \tau_2} B^n(\tau_1,\tau_2) \right]
\]
is given by
$
\hat{y} = \dfrac{n}{\sqrt{2\rho}}.$ \end{corollary}
\begin{remark}
It is natural to ask whether first hitting points are optimal among all two parameter stopping points. More precisely, does the following equality hold?
\[
\sup_{\tau \in \mathcal{T}_y} \mathbb{E}\left[ e^{-\rho \tau_1 \tau_2} h(B(\tau_1, \tau_2)) \right] 
\overset{?}{=}
\sup_{\tau \in \mathcal{T}} \mathbb{E}\left[ e^{-\rho \tau_1 \tau_2} h(B(\tau_1, \tau_2)) \right]?
\]
This remains an open question in our context and will not be addressed in this work.
\end{remark}

\begin{remark}
Our formulation is closely related to the optimal stopping problem studied by Mazziotto \& Szpirglas \cite{mazziotto05}, where the discounted reward functional $$ e^{-\rho t x} h(B(t,x)) $$ is optimized over two parameter stopping sets. Their analysis relies on the theory of Snell envelopes and superharmonic majorants in the plane, leading to a general existence theory for bounded rewards. In contrast, our approach focuses on a subclass of stopping points - first hitting point of levels - and allows unbounded rewards under suitable growth conditions. This yields explicit optimal thresholds and complements the general framework of \cite{mazziotto05} with constructive results.
\end{remark}

\begin{figure}[H]
    \centering
    \includegraphics[width=0.9\linewidth]{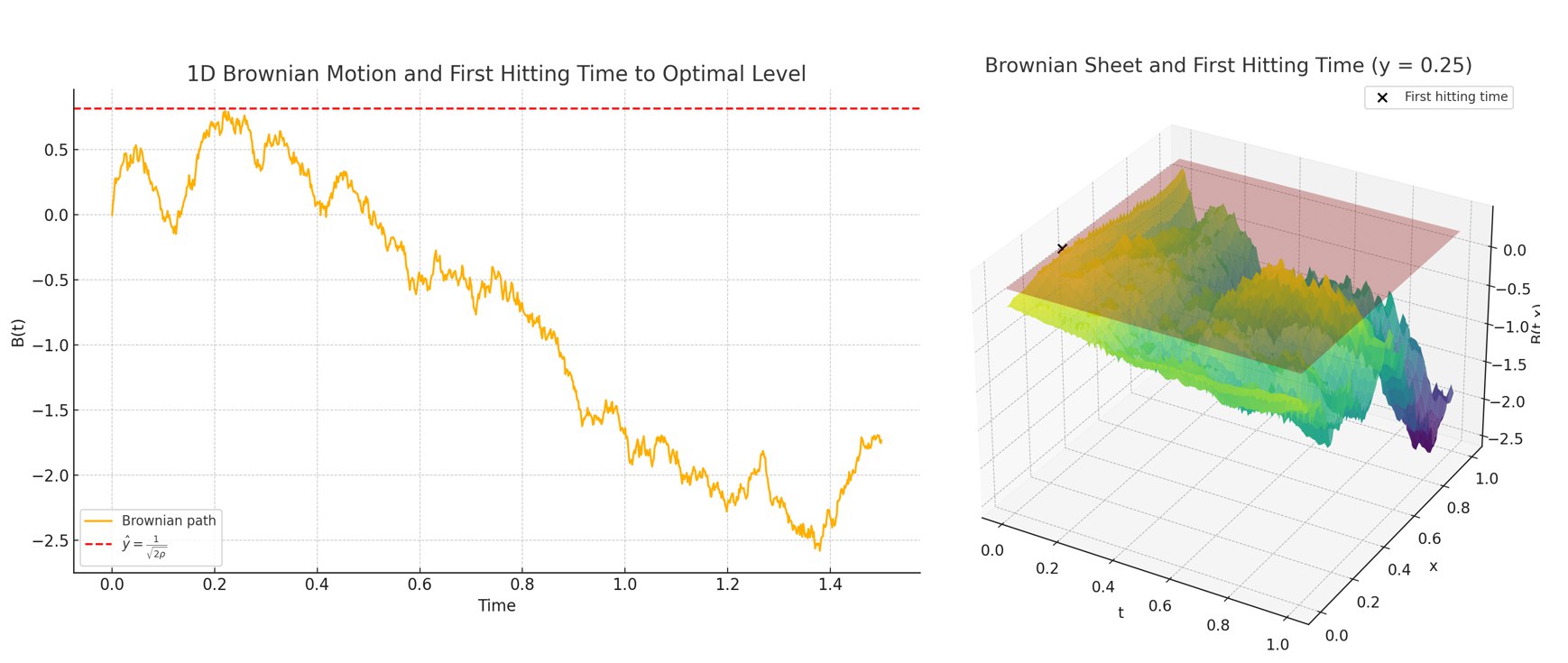}
     \caption{ Comparison of first hitting times/points in the one parameter Brownian motion (left) and the two parameter Brownian sheet (right), both with threshold $y = \hat{y}$. }
    \label{comparison}
\end{figure}

To validate the theoretical findings from the first optimal stopping problem for the Brownian sheet, we perform a Monte Carlo simulation. We generate multiple independent realizations of the Brownian sheet and, for each path, we record the first hitting point \( \tau = (\tau_1, \tau_2) \) of a fixed level \( y = \hat{y} = \frac{1}{\sqrt{2\rho}} \). At this stopping point, we compute the corresponding discounted reward \( e^{-\rho \tau_1 \tau_2} B(\tau_1, \tau_2) \).


\section{The second optimal stopping problem:\\ The integrated discounted Brownian sheet}
In this section we study the two parameter analogue of Problem \ref{2.2}: 
\begin{problem}{(Optimal stopping of the integrated discounted Brownian sheet)}\label{5.1}\\
    Let $\rho >0$ be constant. Find a stopping point $\widehat{\tau}:=(\widehat{\tau_1},\widehat{\tau_2})\in \widehat{\T}$ such that
    \begin{align}
  \sup_{(\tau_1,\tau_2)\in\widehat{\T}} \E[\int_0^{\tau_1}\int_0^{\tau_2} e^{-\rho t x} B(t,x)dt dx ]   =\E[\int_0^{\widehat{\tau}_1}\int_0^{\widehat{\tau}_2} e^{-\rho t x} B(t,x)dt dx ].  
    \end{align}
\end{problem}
The analysis below shows that any maximizer, if it exists, must lie on the positive half-line.

\begin{lemma}\label{lem:concavity}
Let $\rho>0$ and define
\[
F(y)
=
\frac{1}{\rho}
\int_\rho^\infty
\frac{y e^{-y\sqrt{2u}}}{u}\,du,
\qquad y>0.
\]
Then:
\begin{itemize}
\item[(i)] $F\in C^2(0,\infty)$.

\item[(ii)] $\displaystyle \lim_{y\to0^+}F'(y)=+\infty$.

\item[(iii)] $\displaystyle \lim_{y\to\infty}F'(y)=0$ and $F'(y)<0$ for $y$ sufficiently large.

\item[(iv)] $F'$ is strictly decreasing on $(0,\infty)$.
\end{itemize}

Consequently, there exists a unique $y^*>0$ such that
\[
F'(y^*)=0,
\]
and $y^*$ is the unique global maximizer of $F$ on $(0,\infty)$.
\end{lemma}
\dproof
Fix $\rho>0$ and consider
\[
F(y)
=
\frac{1}{\rho}
\int_\rho^\infty
\frac{y e^{-y\sqrt{2u}}}{u}\,du,
\qquad y>0.
\]
Differentiation under the integral sign is justified
by dominated convergence, since the integrand
and its derivatives are exponentially decaying
for $u\ge \rho$.
Hence
\[
F'(y)
=
\frac{1}{\rho}
\int_\rho^\infty
\frac{(1-y\sqrt{2u})e^{-y\sqrt{2u}}}{u}\,du.
\]
As $y\to0^+$,
\[
(1-y\sqrt{2u})e^{-y\sqrt{2u}}
\longrightarrow 1,
\]
and therefore
\[
\lim_{y\to0^+}F'(y)
=
\frac{1}{\rho}
\int_\rho^\infty\frac{1}{u}\,du
=
+\infty.
\]
Thus $F'(y)>0$ for $y$ sufficiently small.\\
As $y\to\infty$, exponential decay dominates
polynomial growth, so
\[
\lim_{y\to\infty}F'(y)=0.
\]
Moreover, the dominant term is negative,
hence $F'(y)<0$ for $y$ sufficiently large.\\
To study the critical points, we use the
equivalent representation
\[
F(y)
=
\frac{2y}{\rho}
E_1\!\bigl(y\sqrt{2\rho}\bigr),
\]
where
\[
E_1(x)=\int_x^\infty\frac{e^{-s}}{s}\,ds.
\]
Then
\[
F'(y)
=
\frac{2}{\rho}
\Bigl(E_1(z)-e^{-z}\Bigr),
\qquad
z=y\sqrt{2\rho}.
\]
Define $g(z)=E_1(z)-e^{-z}$.
Then $g$ is continuous on $(0,\infty)$,
\[
\lim_{z\to0^+}g(z)=+\infty,
\qquad
\lim_{z\to\infty}g(z)=0^-.
\]
Moreover,
\[
g'(z)
=
e^{-z}\Bigl(1-\frac{1}{z}\Bigr),
\]
so $g$ is strictly decreasing on $(0,1)$
and strictly increasing on $(1,\infty)$.
Since $g(1)<0$, the equation $g(z)=0$
admits exactly one solution $z^*>0$.\\
Therefore there exists a unique
\[
y^*=\frac{z^*}{\sqrt{2\rho}}>0
\]
such that $F'(y^*)=0$.
Since $F'$ changes sign from positive
to negative at $y^*$,
this point is the unique global maximizer of $F$.\fproof

\begin{theorem} \label{5.2}
\begin{enumerate}
\item
 For given $y \in \R$ define the first hitting value function by
\begin{align}
\Psi(y) := \sup_{\tau \in \mathcal{T}_y} \E\big[\int_0^{\tau_1}\int_0^{\tau_2} e^{-\rho t x} B(t,x)dt dx \big].
\end{align}
Then\begin{align}
    \Psi(y)= F(y):=\frac{1}{\rho} \int_{\rho}^{\infty} \frac{y e^{-|y| \sqrt{2u}}}{u} du.
\end{align}
Hence $\Psi$ admits a maximum point $y^{*}\in \R$ if and only if $y=y^{*}$ is a maximum point of 
\begin{align}\label{5}
    F(y):= \frac{1}{\rho} \int_{\rho}^{\infty} \frac{y e^{-|y| \sqrt{2u}}}{u} du.
\end{align}
\item 
By Lemma~\ref{lem:concavity}, $F(y)$ has a unique maximizer
$y^*>0$. The maximizer satisfies the first order condition
\[
F'(y^*)=0,
\]
i.e.
 \begin{align}
\rho F'(y^{*})=\int_{\rho}^{\infty} \frac{1}{u} (1-y^{*}\sqrt{2u}) e^{-y^{*}\sqrt{2u}} du = 0.
\end{align}
  
\item
   Let $\widehat{\tau}:=(\widehat{\tau}_1, \widehat{\tau}_2)$ be the first hitting point of $B(t,x)$ for the value $y^{*}$.
     Then $\widehat{\tau}$ is optimal among all first hitting points for Problem \ref{5.1}.
This yields the maximum value
\begin{align}
    \sup_y \Psi(y):= \frac{1}{\rho} \int_{\rho}^{\infty} \frac{y^{*} e^{-y^{*} \sqrt{2u}}}{u} du.
\end{align}
\end{enumerate}
\end{theorem}
  
\begin{remark}
The function $F(y)$ is continuous for $y>0$ and satisfies
\[
\lim_{y\to0^+}F(y)=0,
\qquad
\lim_{y\to\infty}F(y)=0.
\]
Moreover,
\[
F(y)<0 \quad \text{for } y<0,
\qquad
F(y)>0 \quad \text{for } y>0.
\]
Hence any maximizer of $F(y)$ must lie in $(0,\infty)$.
By Lemma~\ref{lem:concavity}, such a maximizer exists and is unique.
\end{remark}

\noindent Proof of Theorem \ref{5.2}:\\
Let $K(t,x):= e^{-\rho t x}$ be the discounting process. Then we see by direct calculation that 
\begin{align}
&\frac{\partial^2}{\partial t \partial x} K(t,x)= -\rho K(t,x)+\rho^2 tx K(t,x),\nonumber\\& \text{ i.e. } \nonumber\\
&dK(z)=(-\rho+\rho^2 tx)K(z) dz; \ z=(t,x).
\end{align}
Therefore, by applying the two dimensional It\^{o} formula (Theorem \ref{itom} in the Appendix) to $f(y_1,y_2)=y_1 y_2$:
\begin{align}
    &K(z) B(z)= \int_{R_z}K(\zeta) B(d\zeta) +\int_{R_z} B(\zeta) dK(\zeta)\nonumber\\
    & +2\int_{R_z}\int_{R_z}(-\rho+\rho^2 sa)K(s,a)ds da B(d\zeta')\nonumber\\
    &= \int_0^t \int_0^x e^{-\rho  s a} B(dsda)+ \int_0^t \int_0^x  B(s,a)(-\rho+ \rho^2sa)e^{-\rho sa} ds da.\nonumber\\
    & +2\int_{R_z}\int_{R_z}(-\rho+\rho^2 sa)e^{-\rho s a}ds da B(d\zeta').
\end{align}
Evaluating at $z=(\tau_1,\tau_2)$, taking expectations, and using that
\[
\E\left[\int_0^{\tau_1}\int_0^{\tau_2} e^{-\rho sa}\,B(ds\,da)\right]=0,
\]
we get,
\begin{align}\label{a}
   \E[ e^{-\rho \tau_1 \tau_2}B(\tau_1,\tau_2)]=\E\big[\int_0^{\tau_1}\int_0^{\tau_2}(-\rho+ \rho^2 sa)e^{-\rho sa} B(s,a) ds da\big]. 
\end{align}
Fix $\tau_1$ and $\tau_2$ and define for $\rho > 0$
\begin{align}
    F(\rho)= \E\big[\int_0^{\tau_1}\int_0^{\tau_2} e^{-\rho s a} B(s,a) ds da\big] \text{ and } G(\rho)= \E\big[ e^{-\rho \tau_1 \tau_2}B(\tau_1,\tau_2)\big].
\end{align}
Then \eqref{a} can be written
\begin{align}
    -\rho^2 F'(\rho) -\rho F(\rho)=  G(\rho),  \text  { or} \nonumber\\
    \rho F'(\rho) + F(\rho) = -\frac{1}{\rho}G(\rho).
\end{align}
This is equivalent to
\begin{align} 
\frac{d}{d\rho}(\rho F(\rho))=  -\frac{1}{\rho}G(\rho),
\end{align}
which gives, for any $\rho_1 > \rho$,
\begin{align}\label{b}
\rho_1 F(\rho_1)= \rho F(\rho)-  \int_{\rho}^{\rho_1} \frac{G(u)}{u} du.
\end{align}
It is easy to see that
\begin{align}
    \lim _{\rho_1 \to \infty} \rho_1 F(\rho_1)=0.
\end{align}
Therefore, letting $\rho_1 \to \infty$ in \eqref{b} we get
\begin{align}
    F(\rho)= \frac{1}{\rho} \int_{\rho}^{\infty} \frac{G(u)}{u} du; \ \rho > 0.
\end{align}
Fix $y\in \R$ and choose $\tau_y=(\tau_1,\tau_2)$ such that $B(\tau_1,\tau_2)=y$ a.s. Then
\[
G(u)=\E\!\left[e^{-u\tau_1\tau_2}B(\tau_1,\tau_2)\right]
= y\,\E\!\left[e^{-u\tau_1\tau_2}\right].
\]
Applying Proposition~\ref{lemma1} with $\beta=\sqrt{2u}$ yields
\[
\E\!\left[e^{-u\tau_1\tau_2}\right]
=\E\!\left[e^{-\frac12(\sqrt{2u})^2\tau_1\tau_2}\right]
=e^{-|y|\sqrt{2u}},
\]
and  $$G(u)=y e^{-|y|\sqrt{2u}}.$$
Hence
\begin{align}
   \E\big[\int_0^{\tau_1}\int_0^{\tau_2} e^{-\rho t x} B(t,x)dt dx \big]= F(\rho)= F(y)= \frac{1}{\rho} \int_{\rho}^{\infty} \frac{y e^{-|y| \sqrt{2u}}}{u} du; \ \rho > 0.
\end{align}
We see that $F(y) < 0$ for all $y < 0$. Therefore the maximum of $y \to F(y)$ over all $y \in \R$, if it exists,  must be attained at some $y^{*} > 0$. If $F(y)$ attains its maximum at $y^* \in \R$, then it satisfies the first order condition
\[
\frac{\partial}{\partial y}F(y^*)=0,
\]
i.e.
\[
\int_\rho^\infty \frac{1}{u}\bigl(1-y^*\sqrt{2u}\bigr)e^{-y^*\sqrt{2u}}du=0.
\]

\fproof

\begin{remark}
We can compute $y^*$ numerically.
With the change of variables $t=\sqrt{2u}$, we obtain
\[
F(y)
=
\frac{1}{\rho}
\int_{\sqrt{2\rho}}^\infty
\frac{2y}{t}
e^{-yt}\,dt,
\qquad y>0.
\]
Thus
\[
F(y)
=
\frac{2y}{\rho}
\int_{\sqrt{2\rho}}^\infty
\frac{e^{-yt}}{t}\,dt.
\]
Introducing $z = y\sqrt{2\rho}$ and writing
\[
\int_{\sqrt{2\rho}}^{\infty}\frac{e^{-yt}}{t}\,dt
=
\int_{z}^{\infty}\frac{e^{-s}}{s}\,ds,
\]
we obtain the representation
\[
F(y)
=
\frac{2y}{\rho}
E_1\!\left(y\sqrt{2\rho}\right),
\]
where
\[
E_1(x)
=
\int_x^\infty \frac{e^{-s}}{s}\,ds
\]
denotes the classical exponential integral.\\
The first-order condition $F'(y^*)=0$ becomes
\[
\int_{\sqrt{2\rho}}^\infty
\frac{1}{t}
(1-yt)e^{-yt}\,dt
=0,
\]
which in terms of $z = y\sqrt{2\rho}$ reduces to
\[
E_1(z) = e^{-z}.
\]
This equation admits a unique positive solution
\[
z^* \approx 0.434818,
\]
and therefore the optimal hitting level is
\[
y^*
=
\frac{z^*}{\sqrt{2\rho}}
\approx
\frac{0.434818}{\sqrt{2\rho}}.
\]
\end{remark}

\begin{figure}[h!]
\centering
\includegraphics[width=0.7\textwidth]{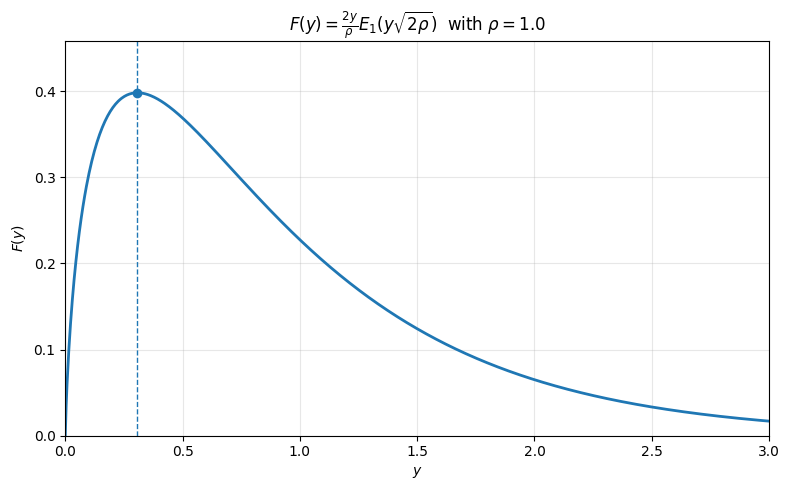}
\caption{Graph of the value function $F(y)$ defined in \eqref{5}.}
\label{fig:Fplot}
\end{figure}

\begin{figure}[t]
\centering

\begin{subfigure}[t]{0.32\textwidth}
\centering
\includegraphics[width=\linewidth]{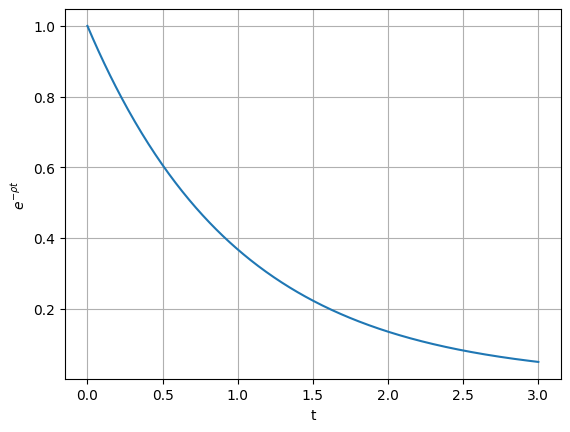}
\caption{One parameter discount $e^{-\rho t}$.}
\label{fig:one-para}
\end{subfigure}
\hfill
\begin{subfigure}[t]{0.32\textwidth}
\centering
\includegraphics[width=\linewidth]{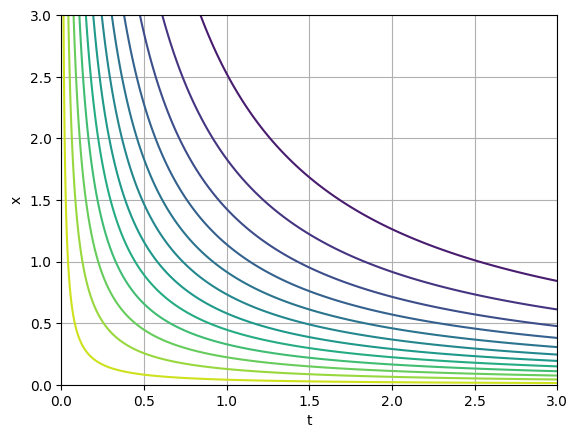}
\caption{Level curves of $e^{-\rho tx}$.}
\label{fig:two-para-contour}
\end{subfigure}
\hfill
\begin{subfigure}[t]{0.32\textwidth}
\centering
\includegraphics[width=\linewidth]{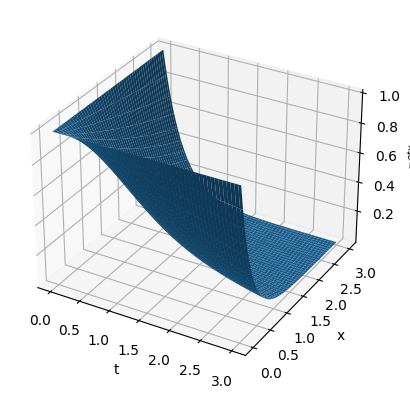}
\caption{Surface of $e^{-\rho tx}$.}
\label{fig:two-para-surface}
\end{subfigure}

\caption{
Geometric comparison of the one and two parameter discount structures.
The multiplicative factor $e^{-\rho tx}$ induces a hyperbolic geometry
and suppresses large rectangles, leading to a qualitatively different
optimal stopping boundary.
}
\label{fig:geometry-comparison}
\end{figure}

\begin{remark}

It is instructive to compare this with the
corresponding one parameter problem (Problem~\ref{2.2}).
In that case the optimal first hitting level is
\[
\widehat y
=
-\frac{1}{\sqrt{2\rho}},
\]
that is, a negative level. 

It is easy to see that in the one-parameter case it is never optimal to stop at a time $\tau$ such that $B(\tau)>0$. Indeed, if $B(\tau)>0$, then by continuity of the Brownian motion there exists $\varepsilon(\omega)>0$ such that $B(t)>0$ for all $t\in[\tau,\tau+\varepsilon]$. Therefore,
continuing the process for a short time yields a strictly larger value of the integral.

In contrast, in the present two parameter
(Brownian sheet) setting we cannot use this argument, because even if $B(\tau_1,\tau_2) > 0$ for a pair of stopping times $\tau_1,\tau_2$ we may still have $B(s,\tau_2)<0$ for some $s < \tau_1$ or $B(\tau_1,a) <0$ for some $a<\tau_2$. and then the gain by continuing the integral over a slightly larger rectangle is not necessarily positive. 
Financially, this may mean that continuation under negative states becomes more costly, since losses are accumulated over a larger region of future cash flows, while the discounting reduces the benefit of waiting.

Thus the optimal stopping level changes sign
when going from one to two parameter problems.
This sign reversal reflects the fundamentally
different geometry of the two cases.
In the one dimensional case the discount acts
along a single time direction, whereas in the
Brownian sheet case the multiplicative factor
$e^{-\rho tx}$ couples the two parameters and
suppresses the contribution of large rectangles.

The geometric difference between the one and two parameter
problems is illustrated in Figure~\ref{fig:geometry-comparison}.
In the one dimensional case (Figure~\ref{fig:one-para})
the discount depends only on time and preserves spatial symmetry.
In contrast, the multiplicative factor $e^{-\rho tx}$
(Figures~\ref{fig:two-para-contour}-\ref{fig:two-para-surface})
induces a hyperbolic geometry and strongly suppresses
contributions from large rectangles.
\end{remark}

\section{Existence of an optimal stopping point in the plane}

In the classical one parameter setting, optimal stopping problems are often addressed using the Snell envelope and the concept of supermartingales. However, in our two parameter setup driven by a Brownian sheet $ B $, the lack of a standard filtration and time direction renders these tools ineffective. Instead, we adopt a potential-theoretic approach based on superharmonic majorants.

Proceeding as in Section 5, we let \( Y_{t_0,x_0,y} \) be as in \eqref{Y} with $p, q$ being zero, and for $t\ge t_0,x\ge x_0$, let $\P^{t_0,x_0,y}$ we denote the probability law of $Y(t,x)$ conditioned on that $Y(t_0,x_0)=y$,  and  write  $\E^{t_0,x_0,y}$ for the expectation with respect to $\P^{t_0,x_0,y}$.
We simplify the notation by setting $z_0=(t_0,x_0), w=(z_0,y)=(t_0,x_0,y) $ .

Define the process $W(z)$ by
\[
W(z)=\begin{bmatrix} z_0+z \\ Y(z)\end{bmatrix},\quad z\geq 0.
\]
Then
\[
{\rm d}W(z)\!=\! \begin{bmatrix} 1\\
\mbox{\quad} \\ \alpha(z)\end{bmatrix}\kern-2pt {\rm d}z\,
+\kern-1pt\begin{bmatrix} 0\\ \mbox{\quad}\\ \beta(z)\end{bmatrix}\kern-2pt B(dz); \\
\quad W(0)\kern-0.5pt=\kern-1.5pt\begin{bmatrix} z_0\\ \mbox{\quad}
\\y\end{bmatrix} =w.
\]

Given a continuous, non-negative reward function \( g \colon \mathbb{R}^3 \to [0, \infty) \), we define the value function by
\[
g^*(w) = \sup_{\tau \in \mathcal{T}} \mathbb{E}^{w}[g(W(\tau_1,\tau_2))].
\]
In our first existence result we show for It\^{o} sheet coefficients $\alpha$ and $\beta$ satisfying the conditions of Theorem \ref{Charact} that the value function is characterized as the least superharmonic majorant \( \widehat{g} \) of \( g \), i.e., \( g^* = \widehat{g} \). Furthermore, the optimal stopping point is given as the first exit point from the continuation region $$ D = \{w : g(w) < \widehat{g}(w)\}.$$ 

This first result is made precise in the following theorem:

\begin{theorem}\label{MainOpt}
Let $g \geq 0$ be a continuous reward function and let $\widehat{g}$ be its least superharmonic majorant. Further, denote by $g^*$ the optimal reward. Then for $\alpha$ and $\beta$ satisfying the conditions of Theorem \ref{Charact} we have:
\begin{itemize}
    \item[(i)] $g^* = \widehat{g}$.
    \item[(ii)] Let $D_\epsilon := \{ w : g(w) < \widehat{g}(w) - \epsilon \}$ for $\epsilon > 0$, and suppose $g$ is bounded. Then for the two parameter exit point $\tau_\epsilon := \tau_{D_\epsilon}$, we have
    \begin{align}\label{10.1.19}
    \left| g^*(w) - \mathbb{E}^{w}[g(W(\tau_{\epsilon}))] \right| \leq \epsilon \quad \text{for all } w.
    \end{align}
    \item[(iii)] In the general case, define the continuation region $D := \{ w : g(w) < g^*(w) \}$. For each $N \geq 1$, define
    \[
    g_N := \min(g, N), \quad D_N := \{ w : g_N(w) < \widehat{g_N}(w) \}, \quad \sigma_N := \tau_{D_N}.
    \]
    Then $D_N \subset D \cap \{ w : g(w) < N \}$ and $D_N \nearrow D$. If $\sigma_N < \infty$ a.e. for all $N$, then
       \[
    g^*(w) = \lim_{N \to \infty} \mathbb{E}^{w}[g(W_{\sigma_N})] \quad \text{for all } w.
    \]
    \item[(iv)] In the special case where $\tau_D < \infty$ a.e. and the sequence $g(W_{\sigma_N})$ is uniformly integrable for each $w$, we have
    \[
    g^*(w) = \mathbb{E}^{w}[g(W_{\tau_D})].
    \]
\end{itemize}
\end{theorem}
The proof will be given in several steps of independent interest:

\subsection{Supermeanvalued, superharmonic and  excessive functions}
We introduce (for general for $\alpha$ and $\beta$) two parameter analogues of the supermartingale property, adapted to the Brownian sheet setting. The first function, which we call a supermeanvalued function, captures the idea that the expected future value of the process, when stopped at a two parameter stopping point, does not exceed its present value.
\begin{definition}
\begin{enumerate}
\item 
A measurable function $f:\mathbb{R}^{3}\to [0,\infty]$ is said to be \emph{supermeanvalued} with respect to $W$ if for all $w=(t_0,x_0,y)$ and all stopping points $\tau_1\geq t_0,\tau_2 \geq x_0$ with $\E[\tau_1 \tau_2] < \infty$ we have
\[
f(w) \geq \mathbb{E}^{w}[f(W(\tau_1,\tau_2))] .
\]
\item
A lower semicontinuous function $f:\mathbb{R}^{3}\to [0,\infty]$ is \emph{excessive} with respect to $W$ if for all $w$ we have
\[
f(w) \geq \mathbb{E}^{w}[f(W(z))] \quad \text{for all } z.
\]
\item
A lower semicontinuous function $f:\mathbb{R}^{3}\to [0,\infty]$ is \emph{superharmonic} with respect to $W$ if for all $w$ and all stopping points $\tau_1\geq t_0,\tau_2 \geq x_0$ with $\E[\tau_1 \tau_2] < \infty$ we have
\[
f(w) \geq \mathbb{E}^{w}[f(W(\tau_1,\tau_2))] .
\]
\end{enumerate}
\end{definition}

We first work in the general state space $w=(t_0,x_0,y)\in\mathbb{R}^3$, so superharmonic and excessive functions are initially defined on $\mathbb{R}^3$.

\begin{lemma}\label{properties}
Let $f, f_1, f_2$, $\{f_j\}_{j\in J}$ be superharmonic or supermeanvalued. Then:
\begin{itemize}
\item[(i)] If $\alpha > 0$, then $\alpha f$ is superharmonic (supermeanvalued).
\item[(ii)] $f_1 + f_2$ is superharmonic (supermeanvalued).
\item[(iii)] $f(w) := \inf_{j \in J} f_j(w)$ is superharmonic (supermeanvalued) for any index set $J$.
\item[(iv)] If $f_k \nearrow f$ pointwise and each $f_k$ is superharmonic (supermeanvalued), then $f$ is superharmonic (supermeanvalued).
\end{itemize}
\end{lemma}

\dproof
See the proof of Lemma 10.1.3 in \cite{O}.
\fproof

\begin{definition}
Let $f \geq h$, where $f$ is superharmonic and $h$ is a real-valued measurable function. Then $f$ is called a \emph{superharmonic majorant} of $h$. The \emph{least superharmonic majorant} is defined by
\[
\overline{h}(w) := \inf \left\{ f(w) : f \text{ is a superharmonic majorant of } h \right\}, \quad w \in \mathbb{R}^3.
\]
\end{definition}

\begin{remark}
Let $g\geq 0$ and let $f$ be a supermeanvalued majorant of $g$. Then if
$\tau$ is a stopping point

$$
f(w)\geq \E^w[f(W(\tau))]\geq \E^w[g(W(\tau))]\;.
$$
So
$$
f(w)\geq\sup_\tau \E^w[g(W(\tau))]=g^\ast(w)\;.
$$
Therefore we always have
\begin{equation} \label{ghat}
\widehat{g}(w)\geq g^\ast(w)\qquad
\text {for all } w.
\end{equation}
What is not so easy to see is that the converse inequality also holds,
i.e.\ that in fact

\begin{equation}
\widehat{g}=g^\ast\;.
\end{equation}
\end{remark}
We will prove this after first establishing some auxiliary results of independent interest:

\begin{theorem}
\label{Charact} Assume that $dW(z)=(dz,\beta (z)B(dz))%
{\acute{}}%
$ for a square integrable predictable process $\beta $ such that $%
(z\longmapsto \mathbb{E}\left[ (\beta (z))^{2}\right] )$ is continuous, $%
\beta (0)=0$ $\mathbb{P}-$a.e. and that  
\begin{equation}
\frac{1}{x}\int_{0}^{x}\mathbb{E}\left[ (\beta (0,v))^{2}\right] dv,x>0
\label{Unbounded}
\end{equation}%
is unbounded. Then a lower semicontinuous function $f:\left[ 0,\infty
\right) \times \mathbb{R}\longrightarrow \left[ 0,\infty \right] $ is
excessive if and only if $f$ is superharmonic.
\end{theorem}

\dproof 
We prove that excessive functions $f$ are superharmonic. Suppose that $f\in
C_{b}((0,\infty )\times \mathbb{R})$ with continuous extensions of all
partial derivatives to $\left[ 0,\infty \right) \times \mathbb{R}$. For
fixed $x$, we can then apply the one dimensional It\^{o}`s formula with
respect to time $t$ and obtain for $Z(t,x):=\int_{0}^{t}\int_{0}^{x}\beta
(u,v)B(du,dv)$ that%
\begin{align*}
f(t_{0}+tx,\,y_{0}+Z(t,x))
&= f(t_{0},y_{0})
 + \int_{0}^{t} x\,f_{t}\!\bigl(t_{0}+ux,\,y_{0}+Z(u,x)\bigr)\,du \\
&\quad + \int_{0}^{t} f_{y}\!\bigl(t_{0}+ux,\,y_{0}+Z(u,x)\bigr)\,Z(du,x) \\
&\quad + \frac12 \int_{0}^{t} f_{yy}\!\bigl(t_{0}+ux,\,y_{0}+Z(u,x)\bigr)
        \left(\int_{0}^{x} (\beta(u,v))^{2}\,dv\right)\,du .
\end{align*}
It follows that
\begin{align}
\E\!\left[f\!\left(t_{0}+tx,\,y_{0}+Z(t,x)\right)\right]
&= f(t_{0},y_{0})  \nonumber\\
&\quad + \E\!\left[
\int_{0}^{t}
\left(
x\,f_{t}\!\left(t_{0}+ux,\,y_{0}+Z(u,x)\right)
\right.\right. \nonumber\\
&\qquad\qquad
\left.\left.
+ \frac12\, f_{yy}\!\left(t_{0}+ux,\,y_{0}+Z(u,x)\right)
\int_{0}^{x} \beta(u,v)^{2}\, dv
\right) du
\right].
\end{align}
Consequently,
\begin{align}
\E\!\left[
\int_{0}^{t}
\left(
x\,f_{t}\!\left(t_{0}+ux,\,y_{0}+Z(u,x)\right)
+ \frac12\, f_{yy}\!\left(t_{0}+ux,\,y_{0}+Z(u,x)\right)
\int_{0}^{x} \beta(u,v)^{2}\, dv
\right) du
\right]
\le 0,
\end{align}
for all $t>0$. Dividing by $t$ and letting $t\downarrow 0$, we obtain
\begin{equation*}
x\,f_{t}(t_{0},y_{0})
+\frac12\,f_{yy}(t_{0},y_{0})\int_{0}^{x}\E\!\left[(\beta(0,v))^{2}\right]\,dv
\le 0
\end{equation*}
for all $t_{0},y_{0}$ and $x$. Dividing by $x>0$ on both sides, we get%
\begin{equation*}
f_{t}(t_{0},y_{0})+\frac{1}{2}f_{yy}(t_{0},y_{0})\frac{1}{x}\int_{0}^{x}%
\mathbb{E}\left[ (\beta (0,v))^{2}\right] dv\leq 0
\end{equation*}%
$t_{0},y_{0}$ and $x>0$. So for $x\searrow 0$ it follows from our
assumptions that%
\begin{equation*}
f_{t}(t_{0},y_{0})\leq 0
\end{equation*}%
for all $t_{0},y_{0}$. Further, the unboundedness of $\frac{1}{x}\int_{0}^{x}%
\mathbb{E}\left[ (\beta (0,v))^{2}\right] dv$, $x>0$ combined with the
boundedness of $f_{t}$ shows that%
\begin{equation*}
f_{yy}(t_{0},y_{0})\leq 0
\end{equation*}%
for all $t_{0},y_{0}$. Thus, we see that $f$ is non-increasing with respect
to $t$ and concave with respect to $y$. Let now $\tau =(\tau _{1},\tau _{2})$
be a bounded stopping point. Then it follows from the
monotonicity property of $f$ with respect to $t$ and Jensen`s inequality for
concave functions that%
\begin{eqnarray*}
&&\mathbb{E}\left[ f(t_{0}+\tau _{1}\tau _{2},y_{0}+Z(\tau _{1},\tau _{2}))%
\right]  \\
&\leq &\mathbb{E}\left[ f(t_{0},y_{0}+Z(\tau _{1},\tau _{2})\right] \leq
f(t_{0},y_{0}+\mathbb{E}\left[ Z(\tau _{1},\tau _{2})\right] ) \\
&=&f(t_{0},y_{0})
\end{eqnarray*}%
for all $t_{0},y_{0}$, which shows that such $f$ in $C_{b}^{2}$ are
superharmonic. The general case can be obtained by using the same arguments
as in Dynkin, 1965 II, p. 5 \cite{DynkinII}.
\fproof

In this paper we also want to prove an existence theorem for optimal stopping points in the case of time-space-homogeneous SDEs in the plane. For this purpose we need the following results: 

\begin{theorem}
\label{Charact2} Assume that $\sigma \in C_{b}^{1}(\mathbb{R})$ such that $%
\sigma $ is bounded away from zero. Let $X_{t,x}^{x},t,x\geq 0$ be the unique
strong solution to%
\begin{equation}
X_{t,x}^{y}=y+\int_{0}^{t}\int_{0}^{x}\sigma (X_{u,v}^{y})B(du,dv)\text{,}%
t,x\geq 0\text{.}  \label{SDEH}
\end{equation}%
Then%
\begin{equation*}
\begin{array}[b]{l}
\begin{array}[b]{l}
\text{(i) If }f\text{ is superharmonic with respect to }X_{t,x_{0}}^{y},t%
\geq 0\text{ for some }x_{0}>0\text{, then} \\ 
\text{ }f\text{ is also superharmonic with respect to }X_{t,x}^{y},t,x\geq 0%
\text{.}%
\end{array}
\\ 
\text{(ii) The least superharmonic majorant }\widehat{g}\text{ of a
non-negative lower } \\ 
\text{semicontinuous function }g\text{ is given by}
\text{ }\widehat{g}(y)=\inf \left\{ f(y):f\text{ concave, }f\geq g\text{ on }%
\mathbb{R}\right\} \text{.}%
\end{array}%
\end{equation*}
\end{theorem}

\dproof
(i) Without loss of generality let $f$ be bounded. We want to show that $f$
is concave. For this purpose let $a<y<b$ and define the $\left\{ \mathcal{F}%
_{t,x_{0}}\right\} _{t\geq 0}-$first exit time%
\begin{equation*}
\tau _{(a,b)}=\inf \left\{ t>0:X_{t,x_{0}}^{y}\notin (a,b)\right\} \text{.}
\end{equation*}%
Then the one dimensional optional sampling theorem implies that%
\begin{equation*}
\mathbb{E}\left[ f(X_{\tau _{(a,b)}\wedge n,x_{0}}^{y})\right] \leq f(y)%
\text{.}
\end{equation*}%
On the other hand, it follows from the $1-$dimensional It\^{o} formula that
\begin{eqnarray*}
(\max (\left\vert a\right\vert ,\left\vert b\right\vert ))^{2} &\geq &%
\mathbb{E}\left[ (X_{\tau _{(a,b)}\wedge n,x_{0}}^{y})^{2}\right] =y^{2}+%
\mathbb{E}\left[ \int_{0}^{\tau _{(a,b)}\wedge n}\int_{0}^{x_{0}}(\sigma
(X_{u,v}^{y}))^{2}dvdu\right]  \\
&\geq &x_{0}\delta ^{2}\mathbb{E}\left[ \tau _{(a,b)}\wedge n\right] 
\end{eqnarray*}%
for a constant $\delta >0$. So Fatou's Lemma yields%
\begin{equation*}
\mathbb{E}\left[ \tau _{(a,b)}\right] <\infty \text{.}
\end{equation*}%
Since $X_{\tau _{(a,b)}\wedge n,x_{0}}^{y}\underset{n\longrightarrow \infty }%
{\longrightarrow }X_{\tau _{(a,b)},x_{0}}^{y}$ a.e. the lower semicontinuity
of $f$ gives%
\begin{equation*}
\lim \inf_{n\longrightarrow \infty }f(X_{\tau _{(a,b)}\wedge
n,x_{0}}^{y})\geq f(X_{\tau _{(a,b)},x_{0}}^{y})\text{.}
\end{equation*}%
Thus%
\begin{equation*}
\mathbb{E}\left[ f(X_{\tau _{(a,b)},x_{0}}^{y})\right] \leq \lim
\inf_{n\longrightarrow \infty }\mathbb{E}\left[ f(X_{\tau _{(a,b)}\wedge
n,x_{0}}^{y})\right] \leq f(y)\text{.}
\end{equation*}%
We have that $X_{\tau _{(a,b)},x_{0}}^{y}\in \left\{ a,b\right\} $ a.e. and
get that%
\begin{equation*}
y=\mathbb{E}\left[ X_{\tau _{(a,b)},x_{0}}^{y}\right] =a\mathbb{P}(X_{\tau
_{(a,b)},x_{0}}^{y}=a)+b\mathbb{P}(X_{\tau _{(a,b)},x_{0}}^{y}=b)\text{.}
\end{equation*}%
So%
\begin{equation*}
\mathbb{P}(X_{\tau _{(a,b)},x_{0}}^{y}=a)=\frac{b-y}{b-a}\,\,\text{ and }\,\,\mathbb{%
P}(X_{\tau _{(a,b)},x_{0}}^{y}=b)=\frac{y-a}{b-a}\text{.}
\end{equation*}%
Therefore we have%
\begin{equation*}
f(y)\geq \mathbb{E}\left[ f(X_{\tau _{(a,b)},x_{0}}^{y})\right] =f(a)\frac{%
b-y}{b-a}+f(b)\frac{y-a}{b-a}\text{,}
\end{equation*}%
which gives 
\begin{equation*}
f((1-\lambda )a+\lambda b)\geq (1-\lambda )f(a)+\lambda f(b)
\end{equation*}%
for $\lambda :=\frac{y-a}{b-a}\in (0,1)$. Hence, $f$ is concave. From
Jensen`s inequality for concave functions (combined with a $2$-parameter
optional sampling theorem for (accessible) stopping points) it then follows
that $f$ is superharmonic with respect to $X_{t,x}^{y},t,x\geq 0$.

(ii) By construction $\widehat{g}$ is the least concave majorant of $g$,
which is lower semicontinuous. So the proof follows from Jensen`s inequality
for concave functions.
\fproof

\begin{corollary}
\label{Corollary}\bigskip Let $X_{t,x}^{y},t,x\geq 0$ be the solution the SDE (%
\ref{SDEH}) and let $g:\mathbb{R}\longrightarrow \left[ 0,\infty \right] $
be a lower semicontinuous function. Denote by $\widehat{g}_{x_{0}}$ the
least superharmonic majorant of $g$ with respect to the process $%
X_{t,x_{0}}^{y},t\geq 0$ and by $\widehat{g}$ the least superharmonic
majorant of $g$ with respect to $X_{t,x}^{y},t,x\geq 0$. Then%
\begin{equation*}
\widehat{g}_{x_{0}}=\widehat{g}\text{.}
\end{equation*}
\end{corollary}

\dproof
We know from Theorem \ref{Charact2} that $\widehat{g}_{x_{0}}$ also is
superharmonic with respect to $X_{t,x}^{y},t,x\geq 0$. So $\widehat{g}\leq 
\widehat{g}_{x_{0}}$. On the other hand, in particular $\widehat{g}$ is
superharmonic with respect to $X_{t,x_{0}}^{y},t\geq 0$. Hence $\widehat{g}%
\geq $ $\widehat{g}_{x_{0}}$.
\fproof

\begin{lemma}
\label{Ex}Let $G$ be an open set in $\mathbb{R}$ and denote by $\tau _{G}$ the
first exit time from $G$ with respect to the one parameter process $%
X_{t,x_{0}}^{y},t\geq 0$ for some fixed $x_{0}$, where $X_{t,x}^{y},t,x\geq 0
$ is the solution of the SDE in Theorem \ref{Charact2}. Assume that $\tau
_{G}<\infty $ a.e. Further, let $h$ be a superharmonic function with respect
to $X_{t,x}^{y},t,x\geq 0$. Then the function $f$ defined by%
\begin{equation*}
f(y)=\mathbb{E}\left[ h(X_{\tau _{G},x_{0}}^{y})\right] 
\end{equation*}%
is superharmonic with respect to $X_{t,x}^{y},t,x\geq 0$.
\end{lemma}

\dproof
Without loss of generality let $h$ be bounded. As in the proof of Theorem %
\ref{Charact2} let $a<y<b$ and define the (integrable) $\left\{ \mathcal{F}%
_{t,x_{0}}\right\} _{t\geq 0}-$first exit time%
\begin{equation*}
\tau _{(a,b)}=\inf \left\{ t>0:X_{t,x_{0}}^{y}\notin (a,b)\right\} \text{.}
\end{equation*}%
Just as in the proof of Theorem \ref{Charact2} we find that%
\begin{equation}
\mathbb{E}\left[ f(X_{\tau _{(a,b)},x_{0}}^{y})\right] =f(a)\frac{b-y}{b-a}%
+f(b)\frac{y-a}{b-a}.  \label{C}
\end{equation}%
On the other hand, we can apply the Theorem of Dambis-Dubins-Schwarz to the $%
\left\{ \mathcal{F}_{t,x_{0}}\right\} _{t\geq 0}-$martingale $%
M_{t}:=X_{t,x_{0}}^{y}-y,t\geq 0$ and get that%
\begin{equation*}
M_{t}=B_{A(t)}^{y}\text{ for all }t\text{ a.e.,}
\end{equation*}%
where $B^{y}$ is a Brownian motion starting in zero (depending on $y$) and%
\begin{equation*}
A(t):=\left\langle M\right\rangle _{t}=\int_{0}^{t}\int_{0}^{x_{0}}(\sigma
(X_{u,v}^{y}))^{2}dudv\text{, }t\geq 0
\end{equation*}%
which increases strictly to infinity a.e. (since $\sigma $ is bounded away
from zero). So%
\begin{equation*}
X_{t,x_{0}}^{y}=y+B_{A(t)}^{y}
\end{equation*}%
Denote by $T$ the continuous inverse of $A$. Then%
\begin{equation*}
X_{T(s),x_{0}}^{y}=y+B_{s}^{y}\text{ for all }s\text{ a.e.}
\end{equation*}%
Let $\tau _{G}^{B}=\tau _{G}^{B}(y)$ be the first exit time from $G$ with
respect to $y+B^{y}$ and consider the first exit time $\tau _{G}$ with
respect to $X_{t,x_{0}}^{y},t\geq 0$ . Let $y\in G$. Then%
\begin{equation*}
X_{\tau _{G},x_{0}}^{y}=y+B_{A(\tau _{G})}^{y}\in \partial G\text{.}
\end{equation*}%
So $A(\tau _{G})\geq \tau _{G}^{B}$. Further, 
\begin{equation*}
X_{T(\tau _{G}^{B}),x_{0}}^{y}=y+B_{\tau _{G}^{B}}^{y}\in \partial G\text{.}
\end{equation*}%
Hence, $T(\tau _{G}^{B})\geq $ $\tau _{G}$, that is $\tau _{G}^{B}\geq
A(\tau _{G})$. Thus, $\tau _{G}^{B}=A(\tau _{G})$ a.e., which implies%
\begin{equation*}
X_{\tau _{G},x_{0}}^{y}=y+B_{\tau _{G}^{B}}^{y}\text{ a.e.}
\end{equation*}%
Since $y+B_{\tau _{G}^{B}}^{y}$ is a measurable functional of $y+B^{y}$, we
can choose any other Brownian motion $W$ (which doesn%
\'{}%
t depend on $y$) and observe that%
\begin{equation*}
y+B_{\tau _{G}^{B}}^{y}\overset{law}{=}y+W_{\tau _{G}^{W}}\text{,}
\end{equation*}%
for $\tau _{G}^{W}=\tau _{G}^{W}(y):=\inf \left\{ t>0:y+W_{t}\notin
G\right\} $. So%
\begin{equation*}
X_{\tau _{G},x_{0}}^{y}\overset{law}{=}y+W_{\tau _{G}^{W}}\text{.}
\end{equation*}%
Further, it follows from the time-homogeneity and the Markov property that%
\begin{equation*}
y+W_{\tau _{G}^{W}}\overset{law}{=}W_{\tau _{G}^{W_{\ast }}}^{\tau
_{(a,b)}^{\ast },y}\text{,}
\end{equation*}%
where $W_{\cdot }^{s,y}$ is a Brownian motion with $W_{s}^{s,y}=y$ a.e, $%
\tau _{(a,b)}^{\ast }:=\inf \left\{ t>0:y+W_{t}\notin (a,b)\right\} $ and%
\begin{equation*}
\tau _{G}^{W_{\ast }}:=\inf \left\{ t>\tau _{(a,b)}^{\ast }:W_{t}^{\tau
_{(a,b)}^{\ast },y}\notin G\right\} \text{.}
\end{equation*}%
On the other hand, one obtains similarly%
\begin{equation*}
X_{\tau _{(a,b)},x_{0}}^{y}\overset{law}{=}y+W_{\tau _{(a,b)}^{\ast }}\text{.%
}
\end{equation*}%
Therefore we see that%
\begin{eqnarray*}
\mathbb{E}\left[ f(X_{\tau _{(a,b)},x_{0}}^{y})\right]  &=&\mathbb{E}\left[
h(W_{\tau _{G}^{W_{\ast }}}^{\tau _{(a,b)}^{\ast },y+W_{\tau _{(a,b)}^{\ast
}}})\right]  
=\mathbb{E}\left[ h(W_{\tau _{G}^{W_{\ast \ast }}}^{y})\right],
\end{eqnarray*}%
where $\tau _{G}^{W_{\ast \ast }}:=\inf \left\{ t>\tau _{(a,b)}^{\ast
}:y+W_{t}\notin G\right\} $. Since $\tau _{G}^{W_{\ast \ast }}\geq \tau
_{G}^{W}$, it follows from Jensen`s inequality for concave functions and the
optional sampling theorem that 
\begin{equation*}
\mathbb{E}\left[ f(X_{\tau _{(a,b)},x_{0}}^{y})\right] =\mathbb{E}\left[
h(W_{\tau _{G}^{W_{\ast \ast }}}^{y})\right] \leq \mathbb{E}\left[ h(W_{\tau
_{G}^{W}}^{y})\right] =f(y)\text{.}
\end{equation*}%
So the latter and (\ref{C}) imply that $f$ must be concave. Then the proof
follows from (i) of Theorem \ref{Charact2}.

\fproof



 



\subsection{Least superharmonic majorants}

In this section, we recall the constructive procedure for obtaining the least superharmonic majorant of a given reward function, following \cite{O}.

\begin{lemma}\label{least}
Consider a lower semicontinuous function $g = g_0 : \mathbb{R}^3 \to [0,\infty)$ and define a sequence of functions $g_n$, for $n \geq 0$, recursively by
\[
g_n(w) = \sup_{z} \mathbb{E}^{w}\left[ g_{n-1}(W(z)) \right],
\]
where 
\[
S_n := \left\{ k2^{-n} : 0 \leq k \leq 4^n \right\}, \quad 
S_n^* := \left\{ k2^{-n} : 0 \leq k \leq 4^n \right\} \cup \left\{ -k2^{-n} : 0 \leq k \leq 4^n \right\}.
\]
Then $g_n \nearrow \widehat{g}$ pointwise, where $\widehat{g}$ is the least superharmonic majorant of $g$, that is,
\[
\widehat{g}(w) = \overline{g}(w) := \inf \left\{ f(w) : f \text{ supermeanvalued majorant of } g \right\}.
\]
\end{lemma}

\dproof
We follow the arguments of Theorem 10.1.7 in \cite{O}. It is clear that the sequence $(g_n)_{n \geq 0}$ is increasing pointwise. Define
\[
\check{g}(w) := \lim_{n \to \infty} g_n(w).
\]
Then, for all $n \geq 1$ and all $w \in S_n \times S_n^*$, we have
\[
\check{g}(w) \geq g_n(w) \geq \mathbb{E}^w [g_{n-1}(W(z))].
\]
As $S_n, S_n^*$ are increasing with $n$, it follows that
\[
\check{g}(w) \geq \liminf_{n \to \infty} \mathbb{E}^w [g_{n-1}(W(z))] \geq \mathbb{E}^w \left[ \liminf_{n \to \infty} g_{n-1}(W(z)) \right] = \mathbb{E}^w [\check{g}(W(z))],
\]
for all $w$ in the dense subset $S := \bigcup_{n \geq 1} S_n \times S_n^*$. 

Since $g_n$ is lower semicontinuous (see Lemma 8.1.4 in \cite{O}), and limits of increasing lower semicontinuous functions are also lower semicontinuous, we conclude that $\check{g}$ is lower semicontinuous.

Now, let $w \in \mathbb{R}^3$ and choose a sequence $w_k \in S$ such that $w_k \to w$. Then, using Fatou's lemma and lower semicontinuity,
\[
\check{g}(w) \geq \liminf_{k \to \infty} \mathbb{E}^w [\check{g}(W(z_k))] \geq \mathbb{E}^w \left[ \liminf_{k \to \infty} \check{g}(W(z_k)) \right] \geq \mathbb{E}^w [\check{g}(W(z))].
\]
Hence, $\check{g}$ is an excessive function, and by Theorem \ref{Charact}, it is also superharmonic.

Finally, let $f$ be any supermeanvalued majorant of $g$. Then by induction, $f(w) \geq g_n(w)$ for all $n$, so $f(w) \geq \check{g}(w)$. Thus, $\check{g}$ is the least supermeanvalued majorant of $g$, i.e., $\check{g} = \overline{g}$.
\fproof

\subsection{Completing the proof of Theorem \ref{MainOpt}}
We now have all we need to complete the proof of Theorem \ref{MainOpt}.

We follow the proof strategy from Theorem 10.1.9 in \cite{O}.

First assume that $g$ is bounded and define
\begin{equation}
\widetilde{g}_\epsilon(w)=\E^{w}[\widehat{g}(W(\tau_\epsilon))]\qquad
\mbox{for $\epsilon>0$}\;.
\end{equation}
Under our assumptions (based on an It\^{o} sheet) , $W$ is additive in the initial value. So using Fubini`s theorem it then follows that $\widetilde{g}_\epsilon$ is supermeanvalued. We
claim that
\begin{equation}\label{claim}
g(w)\leq\widetilde{g}_\epsilon(w)\qquad
\text{for all } w.
\end{equation}
To see this define
\begin{equation}
\beta\colon=\sup_w\{g(w)-
\widetilde{g}_\epsilon(w)\}\;.
\end{equation}
Then $\widetilde{g}_\epsilon(w)+\beta$ is a supermeanvalued majorant of $g$. Hence
\begin{equation}
\widehat{g}(w)<\widetilde{g}_\epsilon(w)+\beta\quad\mbox{for all $w$}\;.
\end{equation}
Choose $\alpha$ such that $0<\alpha<\epsilon$. Then there exists $w_0$ such that
$$
\beta-\alpha<g(w_0)-\widetilde{g}_\epsilon(w_0)\;.
$$
Clearly
$$
0\leq\widehat{g}(w_0)-g(w_0)\leq\widetilde{g}_\epsilon(w)+\beta-g(w_0)<\alpha\;.
$$
Hence
\begin{equation}
w_0\in D_\alpha^C\subseteq D_\epsilon^C
\end{equation}
and therefore
$$
\widetilde{g}_\epsilon(x)=
\E^{w_0}[\widehat{g}(W(\tau_\epsilon))]=
\widehat{g}(w_0)\;.
$$
This gives
\begin{equation}
\beta-\alpha<g(w_0)-\widetilde{g}_\epsilon(w_0)\leq
\widehat{g}(w_0)-\widetilde{g}_\epsilon(w_0)=0\;.
\end{equation}
Letting $\alpha\downarrow 0$, we conclude that
$$
\beta\leq 0\;,
$$
which proves the claim \eqref{claim}.

We conclude that $\widetilde{g}_\epsilon$ is a supermeanvalued majorant of $g$. Therefore
\begin{equation} \label{geps}
\widehat{g}\leq\widetilde{g}_\epsilon=
\E^w[\widehat{g}(W(\tau_\epsilon))]\leq
\E^w[(g+\epsilon)(W(\tau_\epsilon))]\leq
g^\ast+\epsilon
\end{equation}
and since $\epsilon$ was arbitrary, we have by \eqref{ghat}
$$
\widehat{g}=g^\ast\;.
$$
If $g$ is not bounded, let
$$
g_{_N}=\min(N,g)\;,\qquad N=1,2,\ldots
$$
and as before let $\widehat{g_{_N}}$ be the least superharmonic
majorant of $g_N$. Then
$$
g^\ast\geq g_{_N}^\ast=\widehat{g_{_N}} \uparrow h\qquad
\mbox{as $\,N\to\infty\,$, where $\,h\geq\widehat{g}$}
$$
since $h$ is a superharmonic majorant of $g$. Thus
$h=\widehat{g}=g^\ast$ and  this proves \eqref{claim} for general $g$.\\
 From
\eqref{geps} and \eqref{claim} and  we obtain \eqref{10.1.19}.

Finally, to obtain (iii) and (iv) let us again first assume that $g$ is
bounded. Then, since
$$
\tau_\epsilon\uparrow\tau_D\qquad\mbox{as $\,\epsilon\downarrow 0$}
$$
and $\tau_D<\infty$ a.s we have
\begin{equation}
\E^w[g(W(\tau_\epsilon))]\to
\E^w[g(W(\tau_D))]\qquad \mbox{as $\,\epsilon\downarrow 0$}\;,
\end{equation}
and hence by \eqref{geps} and \eqref{claim}
\begin{equation}\label{gbounded}
g^\ast(w)=\E^w[g(W(\tau_D))]\qquad\mbox{if $g$ is bounded}\;.
\end{equation}
Finally, if $g$ is not bounded define
$$
h:=\lim_{N\to\infty}\widehat{g_{_N}}\;.
$$
Then $h$ is superharmonic by Lemma \ref{properties} (iv) and since
$\widehat{g_{_N}}\leq\widehat{g}$ for all $N$ we have
$h\leq\widehat{g}$. On the other hand $g_{_N}\leq \widehat{g_{_N}}\leq
h$ for all $N$ and therefore $g\leq h$. Since $\widehat{g}$ is the
least superharmonic majorant of $g$ we conclude that
\begin{equation} \label{h}
h=\widehat{g}\;.
\end{equation}
Hence by \eqref{gbounded}, \eqref{h} we obtain (iii):
$$
g^\ast(w)=\lim_{N\to\infty}\widehat{g_{_N}}(w)=
\lim_{N\to\infty} \E^w[g_{_N}(W(\sigma_N))]
\leq\lim_{N\to\infty}\E^w[g(W(\sigma_N))]\leq
g^\ast(w)\;.
$$
Note that $\widehat{g_{_N}}\!\leq\! N$ everywhere, so if
$g_{_N}(w)\!<\!\widehat{g_{_N}}(w)$ then $g_{_N}(w)\!<\!N$ and
therefore
$g(w)=g_{_N}(w)<\widehat{g_{_N}}(w)\leq\widehat{g}(w)$ and
$g_{_{N+1}}(w)=g_{_N}(w)<\widehat{g_{_N}}(w)\leq
\widehat{g_{_{N+1}}}(w)$.
Hence $D_N\subset D\cap\{w;g(w)<N\}$ and $D_N\subset D_{N+1}$ for all
$N$. So by \eqref{h} we conclude that $D$ is the increasing union of the
sets $D_N$; $N=1,2,\ldots$ Therefore
$$
\tau_D=\lim_{N\to\infty}\sigma_N\;.
$$
So by (iii) and uniform integrability we have
\begin{eqnarray*}
\widehat{g}(w) &=& \lim_{N\to\infty}\widehat{g_{_N}}(w)=
    \lim_{N\to\infty} \E^w[g_{_N}(W(\sigma_N))]  \\
&=& \E^w[\lim_{N\to\infty}g_{_N}(W(\sigma_N))]
=\E^w[g(W(\tau_D))]\;,
\end{eqnarray*}
and the proof of Theorem \ref{MainOpt} is complete.
\fproof

\bigskip

We are coming now to our main result on (accessible) optimal stopping
points for time-space homogeneous SDEs in the plane:

\begin{theorem}\label{MainOpt2}
Let $g \ge 0$ be a continuous reward function and let $\widehat g$ 
be its least superharmonic majorant with respect to the solution process 
$X_{t,x}^{y}$, $t,x \ge 0$, of \eqref{SDEH} in Theorem~\ref{Charact2}. 
Let $g^*$ denote the optimal reward.
\medskip
\noindent
\textnormal{(i)} 
We have
\[
g^* = \widehat g .
\]
\medskip
\noindent
\textnormal{(ii)} 
Assume that $g$ is bounded and fix $x_0>0$. 
For $\varepsilon>0$, define
\[
D_\varepsilon 
:= \{x : g(x) < \widehat g(x) - \varepsilon \},
\]
and let $\tau_\varepsilon := \tau_{D_\varepsilon}$ 
be the first exit time of $X_{t,x_0}^y$ from $D_\varepsilon$.
Then for all $y$,
\[
\bigl| g^*(y) - \E\!\left[g\!\left(X_{\tau_\varepsilon,x_0}^{y}\right)\right] \bigr|
\le \varepsilon .
\]
\medskip
\noindent
\textnormal{(iii)} 
For a general $g$, define the continuation region
\[
D := \{x : g(x) < g^*(x)\}.
\]
For $N \ge 1$, set
\[
g_N := \min(g,N), 
\qquad 
D_N := \{x : g_N(x) < \widehat{g_N}(x)\},
\]
and define $\sigma_N := (\tau_{D_N}, x_0)$.
Then
\[
D_N \subset D \cap \{x : g(x) < N\},
\qquad
D_N \nearrow D.
\]
If $\tau_{D_N} < \infty$ almost everywhere for all $N$, then
\[
g^*(y)
=
\lim_{N \to \infty}
\E\!\left[g\!\left(X_{\sigma_N}^{y}\right)\right],
\qquad \text{for all } y.
\]
\medskip
\noindent
\textnormal{(iv)} 
If $\tau_D < \infty$ almost everywhere and the sequence 
$g(X_{\sigma_N}^{y})$, $N \ge 1$, is uniformly integrable for all $y$, 
then
\[
g^*(y)
=
\E\!\left[g\!\left(X_{\tau}^{y}\right)\right],
\qquad
\tau := (\tau_D, x_0).
\]
\end{theorem}

\begin{remark}
Note $\tau :=(\tau _{D},x_{0})$ is indeed a strong stopping point, since $%
\tau _{D}$ is a stopping time with respect to the filtration $\mathcal{F}%
_{t,x_{0}},t\geq 0$.
\end{remark}

\dproof
We partially follow the proof of Theorem \ref{MainOpt}. We know from Lemma %
\ref{Ex} for $h:=\widehat{g}$ and $\tau _{G}:=\tau _{\epsilon }=\tau
_{\epsilon }(x_{0})$ that the function $\widetilde{g}_{\epsilon }$ defined by%
\begin{equation*}
\widetilde{g}_{\epsilon }(y)=\mathbb{E}\left[ \widehat{g}(X_{\tau _{\epsilon
},x_{0}}^{y})\right] \text{ for }\epsilon >0
\end{equation*}%
is superharmonic. Further Corollary \ref{Corollary} gives $\widehat{g}%
_{x_{0}}=\widehat{g}$. From here we can proceed further and use the same
arguments as in the proof of Theorem \ref{MainOpt} to obtain the result.
\fproof

\textbf{Acknowledgments}
Nacira Agram and Olena Tymoshenko gratefully acknowledge the financial support provided by the Swedish Research Council grants (2020-04697) and the MSCA4Ukraine grant (AvH ID:1233636), which is funded by the EU, respectively.

\bibliographystyle{plain} 
\bibliography{References}

\section{Appendix: A multidimensional two parameter It\^o formula}
We present a formulation of the multidimensional It\^o formula for
two parameter stochastic integrals adapted to our setting,
building on the work of Wong and Zakai~\cite{WZ}.

\begin{theorem} \label{ito}
[Wong--Zakai type multidimensional It\^o formula]\label{itom} \\
Suppose $Y(z)=Y(t,x) \in \R^n$ is a two parameter stochastic process given as the solution of the following SPDE:
\small
\begin{align}\label{Y3}
Y(z)&=\left[\begin{array}{clcr}
Y_1(z)\\Y_2(z)\\...\\Y_n(z)
\end{array} \right]
=\left[\begin{array}{clcr}
Y_1(0) \\Y_2(0)\\ ... \\Y_n(0)
\end{array} \right] 
+\int_{R_z} \left[ \begin{array}{c}
\alpha_1(\zeta)\\ \alpha_2(\zeta)\\ ...  \\\alpha_n(\zeta)
\end{array} \right] d\zeta
+\int_{R_z}\left[ \begin{array}{rc}
\beta_1(\zeta) \\ \beta_2(\zeta) \\ ... \\ \beta_n(\zeta) 
\end{array} \right]
\left[\begin{array}{rc} B_1(d\zeta)\\ B_2(d\zeta) \\ ... \\B_m(d\zeta)
\end{array} \right],
\end{align}
where
\begin{itemize}
    \item $R_z=[0,t]\times[0,x]$, $z=(t,x)$, $\alpha_k(z) \in \mathbb{R} \text{ for all } k=1,2, ...,n;$  
    \item $\beta_{\ell}(z)= (\beta_{\ell,1}(z), \beta_{\ell,2}(z), ... , \beta_{\ell,m}(z)) \in \mathbb{R}^{1 \times m}$ 
 is an $m$-dimensional row vector for all  $\ell = 1,2, ..., n; $ 
 \item  $B(z)=(B_1(z), B_2(z), ... , B_m(z))^{T} \in \mathbb{R}^{m \times 1}$ is an $m$-dimensional Brownian sheet. 
\end{itemize}

Then, if $f:\mathbb{R}^n \mapsto \mathbb{C}$ is smooth, we have
\small
\begin{align*}
&f(Y(z)) =f(Y(0))+\int_{R_{z}} \sum_{k=1}^n \frac{\partial f}{\partial y_k}(Y(\zeta))\Big[\alpha_k(\zeta
)d\zeta+\beta_k(\zeta)B(d\zeta)\Big]\\
&+\tfrac{1}{2}\int_{R_{z}}\sum_{k,\ell =1}^n \frac{\partial^2 f}{\partial y_k \partial y_{\ell}} (Y(\zeta))\beta_k(\zeta) \beta_{\ell}^{T}(\zeta)d\zeta\\
& +\iint\limits_{R_{z}\times R_{z}}\sum_{k,\ell =1}^n \frac{\partial ^2 f}{\partial y_k \partial y_{\ell}}(Y(\zeta\vee\zeta^{\prime
}))\beta_k(\zeta)B(d\zeta)\beta_{\ell}(\zeta^{\prime})B(d\zeta^{\prime})\\ 
&+\iint\limits_{R_{z}\times R_{z}}\Big\{\sum_{k,\ell =1}^n \frac{\partial ^2 f}{\partial y_k \partial y_{\ell}}(Y(\zeta\vee\zeta^{\prime
}))\beta_k(\zeta^{\prime}) \alpha
_{\ell}(\zeta)\\
&+\tfrac{1}{2}\sum_{k,\ell,p=1}^n \frac{\partial ^{(3)} f}{\partial y_k \partial y_{\ell} \partial y_p}(Y(\zeta\vee\zeta^{\prime}))
)\beta_k(\zeta^{\prime}) \Big[\beta_{\ell}(\zeta)\beta_p^T (\zeta)\Big]\Big\}d\zeta B(d\zeta^{\prime})\\
&+\iint\limits_{R_{z}\times R_{z}}\Big\{\sum_{k,\ell =1}^n \frac{\partial ^2 f}{\partial y_k \partial y_{\ell}}(Y(\zeta\vee\zeta^{\prime
}))\beta_k(\zeta)\alpha_{\ell}(\zeta^{\prime})\\
&+\tfrac{1}{2}\sum_{k,\ell,p=1}^n \frac{\partial ^{(3)} f}{\partial y_k \partial y_{\ell} \partial y_p}(Y(\zeta\vee\zeta^{\prime}))
\beta_k(\zeta) \Big[\beta_{\ell}(\zeta^{\prime})\beta_p^T (\zeta^{\prime})\Big]\Big\}B(d\zeta)d\zeta^{\prime}\\
& +\iint\limits_{R_{z}\times R_{z}}I(\zeta \bar{\wedge} \zeta^{\prime}) \Big\{\sum_{k,\ell =1}^n \frac{\partial ^2 f}{\partial y_k \partial y_{\ell}}(Y(\zeta\vee\zeta^{\prime
}))\alpha_k(\zeta^{\prime}) \alpha_{\ell}(\zeta)\\
&
+\tfrac{1}{2}\sum_{k,\ell,p=1}^n \frac{\partial ^{(3)} f}{\partial y_k \partial y_{\ell} \partial y_p}(Y(\zeta\vee\zeta^{\prime}))
\Big[\alpha_k(\zeta^{\prime}) \beta_{\ell}(\zeta)\beta_p^T (\zeta)+\alpha_k(\zeta) \beta_{\ell}(\zeta^{\prime}) \beta_p^{T}(\zeta^{\prime})\Big] \\
& +\tfrac{1}{4} \sum_{k,\ell,p,q=1}^n \frac{\partial ^4 f}{\partial y_k \partial y_{\ell} \partial y_p \partial y_q}(Y(\zeta\vee\zeta^{\prime}))\beta_k(\zeta^{\prime})\beta_\ell^{T}(\zeta^{\prime}%
)\beta_p(\zeta)\beta_q^{T}(\zeta)\Big\}d\zeta d\zeta^{\prime}.
\end{align*}
\end{theorem}

\end{document}